\documentclass[ preprint]{imsart}

\usepackage{amsmath,amsfonts,amssymb,mathrsfs,amsthm}
\RequirePackage{natbib}

\usepackage{amsthm,amsmath,natbib}
\RequirePackage[dvips]{hyperref}

\startlocaldefs
\endlocaldefs

\begin{document}

\bibliographystyle{acmtrans-ims}

\begin{frontmatter}

\title{A Discrete Construction for Gaussian Markov Processes}
\runtitle{A Discrete Construction for Gaussian Markov Processes}


\author{\fnms{Thibaud} \snm{Taillefumier}\ead[label=e1]{ttaillefum@rockefeller.edu}}
\address{Laboratory of Mathematical Physics\\ 1230 York Avenue, New York, NY 10065, USA \\ \printead{e1}}
\affiliation{The Rockefeller University}

\author{\fnms{Antoine} \snm{Toussaint}\ead[label=e2]{atoussai@stanford.edu}}
\address{Financial Math Program\\ 390 Serra Mall, Stanford, CA 94305, USA \\ \printead{e2}}
\affiliation{Stanford University}

\runauthor{Thibaud Taillefumier}

\begin{abstract}
In the L\'evy  construction of Brownian motion, a Haar-derived basis of  functions is used to form a finite-dimensional process $W^{N}$ and to define the Wiener process  as the almost sure path-wise limit of  $W^{N}$ when $N$ tends to infinity. 
We generalize such a construction to the class of centered Gaussian Markov processes $X$ which can be written $X_{t} = g(t) \cdot \int_{0}^{t} f(t) \, dW_{t}$ with $f$ and $g$ being continuous functions.
We build the finite-dimensional process $X^{N}$ so that it gives an exact representation of the conditional expectation of $X$ with respect to the filtration generated by ${\lbrace X_{k/2^{N}}\rbrace}$ for $0 \! \leq \! k \! \leq \! 2^{N}$.
Moreover, we prove that the process $X^{N}$ converges in distribution toward $X$.
\end{abstract}

\begin{keyword}[class=AMS]
\kwd[Primary ]{60G15}
\kwd{60J25}
\kwd{28C20}
\end{keyword}

\begin{keyword}
\kwd{Gaussian Process}
\kwd{Markov Process}
\kwd{Discrete Representation}
\end{keyword}

\end{frontmatter}


\section{Introduction}

Given some probability space, it is often challenging to establish results about continuous adapted stochastic processes. 
As a matter of fact, the mere existence of such processes can prove lengthy and technical:
the direct approach to build continuous Markov processes consists in evaluating the desired finite-dimensional distributions of the process, and then constructing the measure associated with the process on an appropriate measurable space, so that this measure consistently yields the expected finite-dimensional distributions~\cite{Rogers}.\\
In that respect, it is advantageous to have a discrete construction of a continuous stochastic process. 
For more general purpose, a discrete representation of a continuous process proves very useful as well.
Assuming some mode of probability convergence,  at stake is to write a process $X$ as a convergent series of random functions $f_{n} \cdot \xi_{n}$ 
\begin{equation}
X_{t} = \sum_{n=0}^{\infty} f_{n} (t)\cdot \xi_{n}  = \lim_{N \to \infty} \sum_{n=0}^{N} f_{n} (t)\cdot \xi_{n} \, ,
\nonumber
\end{equation}
where $f_{n}$ is a deterministic function and $\xi_{n}$ is a given random variable.\\
The L\'evy construction of Brownian motion --later referred as Wiener process-- provides us with a first example of discrete representation for a continuous stochastic process.
Noticing the simple form of the probability density of a Brownian bridge, it is based on completing sample paths by interpolation according to the conditional probabilities of the Wiener process~\cite{Levy}.
More especially, the coefficients $\xi_{n}$ are Gaussian independent and the elements $f_{n}$, called Schauder elements, are obtained by time-dependent integration of the Haar elements.
This latter point is of relevance since, for being a Hilbert system, the introduction of the Haar basis greatly simplify the demonstration of the existence of the Wiener process~\cite{Ciesielski}.\\
From another perspective, fundamental among discrete representations is the Karhunuen-Lo\`eve decomposition.
Instead of yielding a convenient construction scheme, it represents a stochastic process by expanding it on a basis of orthogonal functions~\cite{Karhunen, Loeve}.
The definition of the basis elements  $f_{n}$ depends only on the second-order statistics of the considered process and the coefficients $\xi_{n}$  are pairwise uncorrelated random variables.
Incidentally,  such a decomposition is especially suited to study Gaussian processes because the coefficients of the representation then become Gaussian independent.
For these reasons, the Karhunen-Lo\'eve decomposition is of primary importance in exploratory data analysis, leading to methods referred as ``principal component analysis", ``Hotelling transform" ~\cite{Hotelling} or ``proper orthogonal decomposition"~\cite{Lumley} according to the field of application.
In particular, it was directly applied to the study of stationary Gaussian Markov processes in the theory of random noise in radio receivers~\cite{Kac}.
\\

In view of this, we propose a construction of Gaussian Markov processes using a Haar-like basis of functions.
The class of processes we consider is general enough to encompass commonly studied centered Gaussian Markov processes that satisfy minimal properties of continuity.
We stress the connection with the Haar basis because our basis of decomposition is the exact analog of the Haar-derived Schauder functions used in the L\'evy construction of the Wiener process.
As opposed to the Karhunene-Lo\`eve decomposition, our basis is not made of orthogonal functions but the elements are such that the random coefficients $\xi_{n}$ are always independent and Gaussian with law $\mathcal{N}(0,1)$, i.e. with zero mean and unitary variance.\\
The almost sure path-wise convergence of our decomposition toward a well-defined continuous process is quite straightforward. 
Most of the work lies in proving that the candidate process provides us with an exact representation of a Gaussian Markov process and in demonstrating that our decomposition converges in distribution toward this representation.
Validating the decomposition essentially consists in proving the weak convergence of all the finite-dimensional measures induced by our construction on the Wiener space:
it requires the introduction of an auxiliary orthonormal system of functions in view of using the Parseval relation.
To furthermore establish the convergence in distribution of the representation, we only need demonstrating the tightness of this family of induced measures.
\\

The discrete construction we present displays both analytical and numerical interests for further applications. \\
Analytically-wise, even if it does not exhibit the same orthogonal properties as the Karhunene-Lo\`eve decomposition, our representation can prove as advantageous to establish analytical results about Gaussian Markov processes. 
It is especially noticeable when computing quantities such as the characteristic functional of random processes~\cite{Hida,Donsker} as shown in annex.
This is just an example of how, equipped with a discrete representation, one can expect to make demonstration of properties about continuous Gaussian Markov processes more tractable.
If the measures induced by our decomposition on the classical Wiener space converge weakly toward the measure of a Gaussian Markov process,
we put forward that the convergence of our decomposition is almost sure path-wise toward the representation of a Gaussian Markov process. 
This result contrasts with the convergence in mean of the Karhunene-Lo\`eve decomposition.\\
From another point of view, three Haar-like properties make our decomposition particularly suitable for certain numerical computations: 
all basis elements have compact support on an open interval with dyadic rational endpoints; these intervals are nested and become smaller for larger indices of the basis element, and for any dyadic rational, only a finite number of basis elements is nonzero at that number. 
Thus the expansion in our basis, when evaluated at a dyadic rational, terminates in a finite number of steps. 
These properties suggest an exact schema to simulate sample paths of a Gaussian Markov process $X$ in an iterative ``top-down'' fashion. 
Assuming conditional knowledge of a sample path on the dyadic points of $D_{N} = \lbrace k2^{-N} \vert 0 \leq k \leq 2^{N}\rbrace$, one can decide to further the simulation of this sample path at any time $t$ in $D_{N+1}$ by drawing a point according to the conditional law of $X_{t}$ knowing ${\lbrace X_t \rbrace}_{t \in D_{N}} $, which is simply expressed in the framework of our construction.
It can be used to great advantage in numerical computations such as dychotomic search algorithms for first passage times: 
considering a continuous boundary, we shall present elsewhere a fast Monte-Carlo algorithm that simulates sample-paths with increasing accuracy only in time regions where a first passage is likely to occur. \\


\section{Main Result}

Beforehand, we emphasize that the analytical and numerical advantages granted by the use of our decomposition come at the price of generality, being only suited for Gaussian Markov processes with minimal properties of continuity. 
We also remark that if the Karhunen-Lo\`eve decomposition is  widely used in data analysis, our decomposition mainly provides us with a discrete construction scheme for Gaussian Markov processes.

\newtheorem*{prop}{Proposition}
\begin{prop}
Let $X=\lbrace X_{t}, \mathcal{F}_{t}; 0 \! \leq \! t \! \leq \! 1 \rbrace $  be a real adapted process  on some probability space $\left( \Omega, \mathcal{F}, \mathrm{\bold{P}} \right)$  which takes value in the set of real numbers and let $\xi_{n,k}$ with $n  \! \geq  \!0$ and $0 \! \leq \! k \!  < \! 2^{n}$ be Gaussian random variables  of law $\mathcal{N}(0,1)$ .\\
If there exist some non-zero continuous functions $f$ and $g$ such that  
\begin{equation}
X_{t} = g(t) \cdot \int_{0}^{t} f(t) \, dW_{t} \, , \quad with \quad 0 \! \leq \! t \! \leq \! 1 \, ,
\nonumber
\end{equation}
then there exists  a basis of continuous functions $\Psi_{n,k}$ for $n \! \geq  \! 0$ and $0 \! \leq  \! k \! <  \!2^{n}$  such that the random variable 
\begin{equation}
X_{t}^{N} = \sum_{n=0}^{N} \sum_{ \hspace{5pt} 0 \leq k < 2^{n\!-\!1} } \; \Psi_{n,k}(t) \cdot \xi_{n,k}
\nonumber
\end{equation}
follows the same law as the conditional expectation of $X_{t}$ with respect to the filtration generated by ${\lbrace X_{k/2^{N}}\rbrace}$ for $0 \! \leq \! k \! \leq \! 2^{N}$.
The functions $\Psi_{n,k}$ thus defined have support in $S_{n,k}\!= \! \left[ k \! \cdot \!2^{-n\!+\!1}, (k\!+\!1)2^{-n\!+\!1}\right]$ and admit simple analytical expressions in terms of functions $g$ and $f$.\\
Moreover, the path-wise limit $\lim_{N \to \infty} X^{N}$ defines almost surely a continuous process which is an exact representation of $X$ and we have
\begin{equation}
X^{N} \stackrel{\mathcal{D}}{\longrightarrow} X \, ,
\nonumber
\end{equation}
meaning that the finite-dimensional process $X^{N}$ converges in distribution toward $X$ when $N$ tends to infinity.
\end{prop}

\newtheorem*{rem1}{Remark}
\begin{rem1}
The function $f$ can possibly be zero on a negligible set in $[0,1]$ in the previous proposition.
\end{rem1}

To prove this proposition, the paper is organized as follows. 
We first review some background about Gaussian Markov processes $X$ and their Doob representations as $X_{t} = g(t) \cdot \int_{0}^{t} f(t) \, dW_{t}$. 
Then we develop the rationale of our construction by focusing on the conditional expectations of the process $X$ with respect to the filtration generated by ${\lbrace X_{k/2^{N}}\rbrace}$ for $0 \! \leq \! k \! \leq \! 2^{N}$. 
In the fifth section, we propose a basis of expansion to form the finite-dimensional candidate processes $X^{N}$ and justify the limit process $\overline{X}$ as the almost sure path-wise convergent process $\lim_{N \to \infty} X^{N}$.
In the sixth section, we introduce the auxiliary Hilbert system and prove an important intermediate result.
In the last section, we show that the finite-dimensional processes $X^{N}$ converge in distribution toward $X$ and that $\overline{X}$ is an exact representation of $X$ .


\section{Background on Gaussian Markov Processes}


\subsection{Basic Definitions}

We first define the class of Gaussian Markov processes. 
Let us consider on some probability space $\left( \Omega, \mathcal{F}, \mathrm{\bold{P}} \right)$ a real adapted process $X=\lbrace X_{t}, \mathcal{F}_{t}; 0 \! \leq \! t \! < \! \infty \rbrace $ which takes value in the set of real numbers.
We stress that the index $t$ of the random variable $X_{t}$ runs in the continuous set $\mathbb{R^{+}}$.
For a given realization $\omega$ in $\Omega$, the collection of outcomes $t \mapsto X_{t}(\omega)$ is a sample path of the process $X$.
We only consider processes $X$ for which the sample paths $t \mapsto X_{t}(\omega)$ are continuous.
With these definitions, we are in a position to state the two properties characterizing a Gaussian Markov process.
\begin{enumerate}
\item We say that $X$ is a Gaussian process if, for any integers $k$ and positive reals $t_{1} < t_{2}< \cdots <t_{k}$, the random vector $(X_{t_{1}},X_{t_{2}},\cdots,X_{t_{k}})$ has a joint normal distribution.
\vspace{5pt}
\item We say that $X$ is a Markov process if, for any $s,t \! \geq \! 0$ and $\Gamma \in \mathcal{B}\left( \mathbb{R}\right)$, with $\mathcal{B}\left( \mathbb{R}\right)$ the set of real Borelians, 
\begin{equation}
 \mathrm{\bold{P}} \left( X_{t+s} \in \Gamma \, \vert \, {\mathcal{F}}_{s}\right)=  \mathrm{\bold{P}} \left( X_{t+s} \in \Gamma \, \vert \, X_{s} \right) \, ,
 \nonumber
\end{equation}
which states that the conditional probability distribution of future states $X_{t+s}$, given the present state and all past states ${\mathcal{F}}_{s}$, depends only upon the present state $X_{s}$.
\end{enumerate}
A Gaussian Markov process is a stochastic process that satisfies both Gaussian and Markov properties.
\\
The Wiener process and the Ornstein-Uhlenbeck process are two well-known examples of Gaussian Markov processes.
The Wiener process is defined as the only continuous process $W_{t}$ for which $W_{0}=0$ and the increments $W_{t}-W_{s}$ are independent of  ${\mathcal{F}}_{s}$ and normally distributed with law $\mathcal{N}\left(0,t-s\right)$. 
These requirements naturally place the Wiener process in the class of Gaussian Markov process.  
The Ornstein-Ulhenbeck process can be defined as a solution of the stochastic differential equation of the form
\begin{equation} \label{eq:UStoch}
dX_{t} = \alpha X_{t}\, dt + dW_{t} \quad \mathrm{with} \quad \alpha \in \mathbb{R}\, ,
\nonumber
\end{equation}
We designate $U^{\alpha}_{t}$ the Ornstein-Ulhenbeck process of parameter $\alpha$ starting at $0$ for $t=0$ and
we give its integral expression 
\begin{equation} \label{eq:UInt}
U^{\alpha}_{t}= \int_{0}^{t} e^{\alpha \left( u-t \right)}\,dW_{u} \, .
\end{equation}
The process ${U^{\alpha}}$ naturally appears as a Gaussian Markov process as well:
it is Gaussian for integrating independent Gaussian contributions and Markov for being solution of a first-order stochastic differential equation.
It is known that both processes can be described as discrete processes with an appropriate basis of random functions~\cite{Thibaud}.


\subsection{The Doob Representation}

The discrete construction of the Wiener process and the Ornstein-Uhlenbeck process is likely to be generalized to a wider class of Gaussian Markov processes because any element of this class can be represented in terms of the Wiener process.
By Doob's theorem~\cite{Doob1}, for any Gaussian Markov process $X$, there exist a real non-zero function $g$ and a real function $f$ in $L^{2}_{loc}(\mathbb{R}^{+})$ such that we have the integral representation of $X$
\begin{equation} \label{eq:IntGM}
X_{t} = g(t) \cdot \int_{0}^{t} f(t) \, dW_{t} \, .
\end{equation}
where $W$ is the standard Wiener process.   
If we introduce the non-decreasing function $h$ defined as 
\begin{equation}
h(t) = \int_{0}^{t} f^{2}(u) \, du \, ,
\nonumber
\end{equation}
then, for any $t,s \geq 0$, the covariance of $X$ can be expressed in terms of functions $h$ and $g$
\begin{equation} \label{eq:XCov}
E \left( X_{t} \cdot X_{s} \right) = g(t)g(s)\cdot h\big(\min(t,s)\big) \, .
\end{equation}
The Doob's representation \eqref{eq:IntGM} indicates  that $X$ is obtained from $W$ by a change of variable in time $t \mapsto h(t)$ and, at any time $t$, by a change of variable in space by a time-dependent  factor $x \mapsto g(t)\cdot x$.
The couple of functions $(f,g)$ that intervenes in the Doob's representation of $X$ is not determined univocally.
Yet, one can defined a canonical class of functions $(f,g)$ which are uniquely defined almost surely in $L^{2}_{loc}(\mathbb{R}^{+})$ if we omit their signs~\cite{Hida}. 
Incidentally, we can compare the integral formulation \eqref{eq:IntGM} with expression \eqref{eq:UInt} of the Ornstein-Uhlenbeck process:
we remark that the representation of this Gaussian Markov process is provided by setting $g(t)=e^{\alpha t}$ and $f(t)=e^{-\alpha t}$, which happens to be its canonical representation.


\subsection{Analytical Results} \label{AnalyticalResults}

The discrete construction of Gaussian Markov processes will rely on two analytical results that we detail in the following.\\
First, the Doob's representation allows us to give an analytical expression for the transition kernel $p( X_{t} \scriptstyle = \displaystyle \! x \, \vert X_{t_{0}} \scriptstyle = \displaystyle \! x_{0} )$ of a general Gaussian Markov process.
As the Doob's representation  is a simple change of variables, it is easy to transform the expression of the Wiener transition kernel $p \left(W_{t}\scriptstyle = \displaystyle \! x\,\vert\,W_{t_{0}} \scriptstyle = \displaystyle \! x_{0}\right)$ to establish
\begin{equation} \label{eq:CondProb}
p( X_{t} \scriptstyle = \displaystyle \! x \, \vert X_{t_{0}} \scriptstyle = \displaystyle \! x_{0} )
 = 
 \frac{1}{g(t) \sqrt{2 \pi  \big( h(t) - h(t_{0}) \big)}} 
\cdot
 \exp{ \left( - \frac{  {\left( \frac{x}{g(t)} - \frac{x_{0}}{g(t_{0})} \right)}^{2} }{2  \big( h(t) - h(t_{0} )\big)}  \right)} \,.
 \nonumber
\end{equation}
We have to mention that this expression is only valid if $h(t) \ne h(t_{0})$, otherwise $X$ is deterministic and $p( X_{t} \scriptstyle = \displaystyle \!x \, \vert X_{t_{0}} \scriptstyle = \displaystyle \! x_{0} ) = \delta_{x_{0}}(x)$. \\
Second, we can use this result to evaluate $q( X_{t_{y}} \scriptstyle = \displaystyle \! y \, \vert X_{t_{x}} \! \scriptstyle = \displaystyle \! x, X_{t_{z}} \! \scriptstyle = \displaystyle \! z)$ with $t_{x}\!<\!t_{y}\!<\!t_{z}$, the probability density of $X_{t}$ knowing its values $x$ and  $z$ at two framing times $t_{x}$ and $t_{z}$. 
Because $X$ is a Markov process, a sample path $t \mapsto X_{t}(\omega)$ which originates from $x$ and joins $z$ through $y$ is just the junction of two independent paths: a path originating in $x$ going to $y$ and a path originating from $y$ going to $z$. 
Therefore, after normalization by the absolute probability for a path to go from $x$ to $y$, we have the probability density 
\begin{eqnarray} \label{eq:QuotProb}
q( X_{t_{y}} \scriptstyle = \displaystyle \!y \, \vert X_{t_{x}} \! \scriptstyle = \displaystyle \! x, X_{t_{z}} \! \scriptstyle = \displaystyle \! z) \quad 
= \quad
\frac{ p( X_{t_{y}} \! \scriptstyle = \displaystyle \! y\,  \vert X_{t_{x}} \! \scriptstyle = \displaystyle \! x) \cdot p( X_{t_{z}} \! \scriptstyle = \displaystyle \! z \, \vert X_{t_{y}} \! \scriptstyle = \displaystyle \! y) }{ p( X_{t_{z}} \! \scriptstyle = \displaystyle \! z \, \vert X_{t_{x}} \! \scriptstyle = \displaystyle \! x ) } \,.
\nonumber
\end{eqnarray}
Thanks to the previous expression, we can compute the distribution of $X_{t_{y}}$ knowing $X_{t_{x}}$ and $X_{t_{z}}$, which is expected to be a normal law because we only consider Gaussian processes.
For a general Gaussian Markov process $X$, we refer to that probability law as $\mathcal{N}({}_{ \scriptscriptstyle X \displaystyle} \mu(t_{y}) ,{}_{ \scriptscriptstyle X \displaystyle} \sigma(t_{y})^{2} )$, with mean value ${}_{ \scriptscriptstyle X \displaystyle} \mu(t_{y})$ and variance ${}_{ \scriptscriptstyle X \displaystyle} \sigma(t_{y})^{2} $.
We show in annex that these parameters satisfy
\begin{equation} \label{eq:muX}
{}_{ \scriptscriptstyle X \displaystyle} \mu(t_{y}) 
= 
\frac{g(t_{y})}{g(t_{x})} \cdot \frac{h(t_{z})-h(t_{y}) }{  h(t_{z})-h(t_{x}) } \cdot x 
+ 
\frac{g(t_{y})}{g(t_{z})} \cdot \frac{ h(t_{y})-h(t_{x}) }{ h(t_{z})-h(t_{x}) } \cdot z \, ,
\end{equation}
\begin{equation} \label{eq:sigmaX}
{}_{ \scriptscriptstyle X \displaystyle} \sigma(t_{y})^{2} 
= 
g^{2}(t_{y}) \cdot \frac{ \big( h(t_{y})-h(t_{x}) \big)  \big(h(t_{z})-h(t_{y})  \big) }{ h(t_{z})-h(t_{x}) } \, .
\end{equation}
Once more, we have to mention that these expressions are only valid if $h(t_{y}) \ne h(t_{x})$.
Wether considering a Wiener process or an Ornstein-Uhlenbeck process, the evaluation of \eqref{eq:muX} and \eqref{eq:sigmaX} with the corresponding expression of $g$ and $h$ leads to the already known results~\cite{Thibaud}.


\section{The Rationale of the Construction}


\subsection{Form of the Discrete Representation}

Form now on, we will suppose that the zeros of the function $f$ pertain to a negligible ensemble in $[0,1]$, causing the function $h$ to be strictly increasing.
We will further restrain ourselves to Gaussian Markov processes for which the functions $f$ and $g$ belong to the set of continuous functions on $[0,1]$ denoted $C(0,1)$.
We remark that in such a case, the functions $t \mapsto {}_{ \scriptscriptstyle X \displaystyle} \mu(t)$ and $t \mapsto {}_{ \scriptscriptstyle X \displaystyle} \sigma(t)$ are continuous on $[0,1]$.\\
Bearing in mind the example of the L\'evy construction for the Wiener process, we want to define a basis of continuous functions $\Psi_{n,k}$ in $C(0,1)$ with $0 \! \leq \! k \! < \! 2^{n}$ to form the discrete process
\begin{equation} 
X_{t}^{N} = \sum_{n=0}^{N} \sum_{ \hspace{5pt} 0 \leq k < 2^{n\!-\!1} } \; \Psi_{n,k}(t) \cdot \xi_{n,k} \, ,
\nonumber
\end{equation}
where $\xi_{n,k}$ are independent Gaussian random variables of standard normal law $\mathcal{N}(0,1)$.
We want to chose $\Psi_{n,k}$ so that $t \mapsto X^{N}(\omega)$ converges almost surely toward $t \mapsto X(\omega)$ when $N$ tends to infinity.
Given the continuous nature of the processes  $X_{N}$, we require that the convergence is uniform and normal on $[0,1]$ to ensure the definition of a continuous limit process $\bar{X} = \lim_{N \to \infty} X^{N}$ on $\left[ 0,1\right]$. 
Moreover, we want to define $\Psi_{n,k}$ on supports $S_{n,k}$ of the form 
\begin{equation}
S_{n,k}\!= \! \left[ k \! \cdot \!2^{-n\!+\!1}, (k\!+\!1)2^{-n\!+\!1}\right]
\nonumber
\end{equation}
As a consequence, the basis of functions $\Psi_{n,k}$ will have the following properties : 
all basis elements have compact support on an open interval with dyadic endpoints; these intervals are nested and becomes smaller for larger indices $n$ of the basis element, and for any dyadic rational, only a finite number of basis elements is nonzero at that number.


\subsection{Conditional Averages of the Process}

Now remains to propose an analytical expression for $\Psi_{n,k}$. 
If we denote $D_{N}$ the set of reals $\lbrace k2^{-N} \, \vert \, 0  \! \leq \! k \! \leq \! 2^{N}\rbrace$, the key point is to consider $Z_{t}^{N} = E(X_{t} \vert {\lbrace X_{s} \rbrace}, s \in D_{ \! N})$ the conditional expectation of the random variable $X_{t}$ given $X_{s}$ with $s$ pertaining to the set of dyadic points $D_{N}$.
The collection of random variables $Z_{t}^{N}$ defined on $\Omega$ specify a continuous random process $Z^{N}$ on $\Omega$.
We notice that, if  $t_{x}\!=\!k2^{-N}$ and  $t_{z}\!=\!(k\!+\!1)2^{-N}$ with $0 \! \leq \! k \! < \! 2^{-N}$ are the two successive points of $D_{N}$ framing t, the random variable $Z_{t}^{N}$ is only conditioned by $X_{t_{x}}$ and $X_{t_{z}}$:
\begin{equation}
Z^{N}_{t} = E(X_{t} \vert {\lbrace X_{s} \rbrace}, s \in D_{ \! N}) = E(X_{t} \vert X_{t_{x}}, X_{t_{z}}) \, .
\nonumber
\end{equation}
Using expression \eqref{eq:muX}, we can express the sample paths $t \mapsto Z^{N}_{t}(\omega)$ as a function of $t$ on $[t_{x},t_{z}]$:
for a given $\omega$ in the sample space $\Omega$, we write
\begin{equation}
Z^{N}_{t}(\omega) = {}_{ \scriptscriptstyle X \displaystyle} \mu_{t_{x},t_{z}}(t, x, z) \stackrel{def}{=} {}_{ \scriptscriptstyle X \displaystyle} \mu^{N,k}(t) \, ,
\nonumber
\end{equation}
where the conditional dependency upon parameters $X_{t_{x}}(\omega)  \! = \! x$ and $X_{t_{z}}(\omega) \! = \! z$ is implicit in ${}_{ \scriptscriptstyle X \displaystyle} \mu^{N,k}$.
The random process  $Z^{N}_{t}$ appears then as a parametric function of ${\lbrace X_{s} \rbrace}_{s \in D_{\!N}}$:
for any $\omega$ in the sample space $\Omega$, the sample path $t \mapsto X_{t}(\omega) $ determines a set of value ${\lbrace x_{s} \rbrace}_{s \in D_{\!N}} = {\lbrace X_{s}(\omega) \rbrace}_{s \in D_{\!N}}$ and by extension a sample path $t \mapsto Z^{N}_{t}(\omega) $ for the process $Z^{N}$.\\
Now two points are worth noticing: 
first, $t \mapsto X_{t}(\omega) $ and $t \mapsto Z_{t}^{N}(\omega)$ are continuous sample paths that coincide on the set $D_{\!N}$; second, we have $\lbrace 0,1\rbrace = D_{0} \subset D_{1} \subset \cdots \subset D_{N}$  a growing sequence of sets with limit ensemble $\mathcal{D}$ the set of dyadic points in $\left[ 0, 1 \right]$, which is dense in $\lbrack 0,1\rbrack$.
Then, provided the path-wise convergence is almost surely normal and uniform on $[0,1]$, the limit process of $Z^{N}$ when $N$ tends to infinity should be continuous and the processes $\lim_{N \to \infty}Z^{N}$ and $X$ should be indistinguishable on $\Omega$.


\subsection{Identification of Conditional Averages and Partial Sums}

Identifying the process $Z^{N}$ with the partial sums $X^{N}$ provides us with a rationale to build the functions $\Psi_{n,k}$.\\
We first need to consider the random variable $Z^{N\!+\!1}_{t}$ on  the support $S_{N\!+\!1,k} \! = \! \left[ t_{x}, t_{z}\right]$ of the function $\Psi_{N\!+\!1,k}$.
The Markov property of the process $X$ entails
\begin{eqnarray}\label{eq:muN+1}
Z^{N\!+\!1}_{t}  &=&  E(X_{t} \vert X_{t_{x}}, X_{t_{y}}, X_{t_{z}}) \nonumber\\
& = & 
\left\{
\begin{array}{ll}
 E(X_{t} \vert X_{t_{x}}, X_{t_{y}}) \quad \textrm{if  \: $t_{x} \leq t  \leq  t_{y}$} \: , \\
 E(X_{t} \vert X_{t_{y}}, X_{t_{z}}) \quad \textrm{if \: $t_{y} \leq t  \leq  t_{z}$} \: .\\
\end{array} \right.
\nonumber
\end{eqnarray}
Hence, for a given $\omega$ the estimation of the sample paths $t \mapsto Z_{t}^{N\!+\!1}(\omega)$ is now dependent upon $y \! = \! X_{t_{y}}(\omega)$ with
$t_{y}$ the midpoint of $t_{x}$ and $t_{z}$\begin{equation}
Z^{N\!+\!1}_{t}(\omega) =   
\left\{
\begin{array}{ll}
{}_{ \scriptscriptstyle X \displaystyle} \mu_{t_{x},t_{y}}(t, x, y) \quad \textrm{if  \: $t_{x} \leq t  \leq  t_{y}$}  \\
{}_{ \scriptscriptstyle X \displaystyle} \mu_{t_{y},t_{z}}(t, y, z) \quad \textrm{if \: $t_{y} \leq t  \leq  t_{z}$} \\
\end{array} \right.
\stackrel{def}{=}  {}_{ \scriptscriptstyle X \displaystyle} \nu^{N,k}(t , y) \, .
\nonumber
\end{equation}
We can identify the conditional process $Z^{N}$ and the partial sums $X^{N}$.
Then, for any $\omega$ in $\Omega$, writing the sample path $t \mapsto Z^{N\!+\!1}_{t}(\omega)$ as a function of $y$, we have
\begin{displaymath}
\begin{array}{cccccc}
\Psi_{N\!+\!1,k}(t) \cdot \xi_{N\!+\!1,k}(\omega) & = &  X_{t}^{N\!+\!1}(\omega) & - & X_{t}^{N}(\omega) & \\
 & = &  Z_{t}^{N\!+\!1}(\omega) & - & Z_{t}^{N}(\omega)  & = {}_{ \scriptscriptstyle X \displaystyle} \nu^{N,k}(t ,y) - {}_{ \scriptscriptstyle X \displaystyle} \mu^{N,k}(t)  \, .
\end{array}
\end{displaymath}
Assuming conditional knowledge on $D_{\!N}$, the quantity $Z^{N}\!=\!E(X_{t} \vert {\lbrace X_{s} \rbrace}, s \in D_{ \! N})$ becomes deterministic and the outcome of the random variable $Z^{N\!+\!1}$ is only dependent upon the values of the process on $D_{\!N\!+\!1} \setminus D_{\!N}$.
More precisely, on the support $S_{N\!+\!1,k}$,  the outcome of $Z_{t}^{N\!+\!1}$ is determined  through the function ${}_{ \scriptscriptstyle X \displaystyle} \nu^{N,k}$ by $y$ the outcome of $X_{t_y}$ given $X_{t_x}\!\!=\!x$ and $X_{t_z}\!\!=\!z$.\\
The distribution of $X_{t_y}$ given $X_{t_x}\!\!=\!x$ and $X_{t_z}\!\!=\!z$ follows the law  $\mathcal{N}({}_{ \scriptscriptstyle X \displaystyle}\mu(t_{y}),$ ${}_{ \scriptscriptstyle X \displaystyle}\sigma(t_{y}))$ and we denote $Y_{N,k}$ a Gaussian variable distributed according to such a law.
With this notation, we are in a position to propose the following criterion to compute the function $\Psi_{N\!+\!1,k}$:
the element $\Psi_{N\!+\!1,k}$ is the only positive function with support included in $S_{N\!+\!1,k}$ such that the random variable $\Psi_{N\!+\!1,k}(t) \cdot \xi_{N\!+\!1,k}$ has the same law as ${}_{ \scriptscriptstyle X \displaystyle} \nu^{N,k}(t ,Y_{N,k}) - {}_{ \scriptscriptstyle X \displaystyle} \mu^{N,k}(t)$.
Direct calculations confirms that the previous relation provides us with a consistent paradigm to define the functions $\Psi_{N\!+\!1,k}$.
Incidentally, we have an interpretation for the statistical contribution of the components $\Psi_{N\!+\!1,k} \cdot \xi_{N\!+\!1,k}$:
if one has previous knowledge of $X_{t}$ on $D_{\!N}$, the function $\sum_{k} \Psi_{N\!+\!1,k} \cdot \xi_{N\!+\!1,k}$ represents the uncertainty about $X_{t}$ that is discarded by the knowledge of its value on $D_{\!N\!+\!1}\setminus D_{\!N}$. 



\section{The Candidate Discrete Process} \label{sec:NormConv}

\subsection{The Basis of Functions}

We recall that we carry out the case for which the function $f$ has a negligible set of zeros in $[0,1]$, which directly follows from the previous section.
Before specifying the candidate basis elements $\Psi_{n,k}$, we introduce the following short notations to simplify the writing of their expressions
\begin{equation}
l_{n,k} = \left(2k\right)2^{-n} \quad , \quad m_{n,k} =  \left(2k\!+\!1\right)2^{-n}  \quad , \quad r_{n,k} = 2\left(k\!+\!1\right)2^{-n} \, .
\nonumber
\end{equation}
Then for $n>0$ and $0 \leq k < 2^{n\!-\!1}$, the explicit formulation of the basis of functions $\Psi_{n,k}$ reads
\begin{equation} \label{eq:PsiDef}
\Psi_{n,k}(t)= \left\{ \begin{array}{lllll}
\displaystyle L_{n,k} \cdot g(t)\big( h(t) - h(l_{n,k})\big)
\quad \textrm{if \: $l_{n,k}  \leq  t  < m_{n,k}$} \: , \\
\\
\displaystyle R_{n,k}  \cdot g(t)\big( h(r_{n,k}) - h(t)\big)
\quad \textrm{if \: $m_{n,k}  \leq t < r_{n,k}$} \: ,\\
\\
 0  \quad \textrm{otherwise} \: , \end{array} \right.
\end{equation}
where we use the constants $L_{n,k}$ and $R_{n,k}$ that are defined by the relations
\begin{equation}
L_{n,k}
=
\sqrt{
\frac{ h(r_{n,k})-h(m_{n,k})}
{ \big( h(r_{n,k})-h(l_{n,k}) \big) 
  \big( h(m_{n,k})-h(l_{n,k}) \big)}
 }
 \nonumber
\end{equation}
\begin{equation}
R_{n,k}
=
\sqrt{
\frac{ h(m_{n,k})-h(l_{n,k})}
{ \big( h(r_{n,k})-h(l_{n,k}) \big) 
  \big( h(r_{n,k})-h(m_{n,k}) \big)}
 }\nonumber
\end{equation}
For $N=0$, the basis element $\Psi_{0,0}$ needs to satisfy the relation
\begin{equation}
\Psi_{0,0}(t) \cdot \xi_{0,0} = E(X_{t} \vert \lbrace X_{s} \rbrace , s \in D_{0}=\{0,1\}) \, ,
\nonumber
\end{equation}
which completely defines the analytical expression of $\Psi_{0,0}$ as follows
\begin{equation}  
\Psi_{0,0}(t) = \frac{g(t) \cdot \big( h(t) - h(l_{0,0})\big)}{\sqrt{h(r_{0,0})-h(l_{0,0})}} \,.
\nonumber
\end{equation}
As expected, we directly ascertain the continuity of the $\Psi_{n,k}$ by continuity of $f$ and $g$.\\
We should briefly discuss the form of the functions $\Psi_{n,k}$.
In the case for which $g(t) = 1$ and $h(t) = t$, we find the usual expression of $\Psi_{n,k}$ for the L\'evy construction of a Wiener process:
the elements of the basis are the triangular wedged-functions obtained by integration of  $H_{n,k}$ the standard Haar functions.
In the general case of a Gaussian Markov process, the expression of $\Psi_{n,k}$ can be derived from the Wiener process basis elements by three operations:
a change of variable in time $dt \mapsto h'(t)dt = f^{2}(t)dt$, a time-dependent change of variable in space $x \mapsto g(t) \cdot x$ and a multiplication by the coefficients $L_{n,k}$ and $R_{n,k}$. 
The effect of this multiplication by $L_{n,k}$ and $R_{n,k}$ will be explain in section \ref{sec:CovCal}.\\
Moreover, the paradigm of the construction makes no assumption about the form of the binary tree of nested compact supports $S_{n,k}$. 
Let us consider a given segment $I_{0,0} = [l_{0,0},r_{0,0}[$ and construct by recurrence such a tree.
We suppose that we have the following partition
\begin{equation}
I_{0,0} = \bigcup_{0 \leq k < 2^{n\!-\!1}} I_{n,k} = \bigcup_{0 \leq k < 2^{n\!-\!1}} \big[l_{n,k},r_{n,k}\big[ \, .
\nonumber
\end{equation}
For each $k$ such that $0 \leq k < 2^{n\!-\!1}$, we draw a point $m_{n,k}$ in $S_{n,k}$. 
Then, we have
\begin{equation}
I_{n,k} = \big[l_{n,k},m_{n,k}\big[ \: \bigcup \: \big[m_{n,k},r_{n,k}\big[ = I_{n\!+\!1,2k} \: \bigcup \: I_{n\!+\!1,2k\!+\!1} \,
\nonumber
\end{equation}
and by construction, we posit $m_{n,k} = l_{n\!+\!1,2k\!+\!1} = r_{n\!+\!1,2k}$. 
Iterating the process for increasing $n$, we build a tree of nested compact supports $S_{n,k}=\overline{I_{n,k}}$. \\
The definition  \eqref{eq:PsiDef} enables us to explicit elements $\Psi_{n,k}$ that are adapted to any such tree.
The so-defined functions $\Psi_{n,k}$ will appear to be valid basis elements to build a discrete representation of $X$ under the only requirement that
\begin{equation} \label{eq:DenseCond}
\lim_{n \to \infty} \; \sup_{0 \leq k < 2^{n\!-\!1}} \vert  r_{n,k} - l_{n,k} \vert = 0 \, .
\end{equation}


\subsection{The almost sure Normal and Uniform Convergence} 

We want to prove the validity of the discrete representation of a Gaussian Markov process with Doob's representation \eqref{eq:IntGM} using the proposed basis of functions $\Psi_{n,k}$. 
Let us consider the partial sums $X^{N}$ defined on $\Omega$ by
\begin{equation} \label{eq:SimpSum}
X^{N}_{t}  = \sum_{n=0}^{N} \sum_{ \hspace{5pt} 0 \leq k < 2^{n\!-\!1} } \; \Psi_{n,k}(t) \cdot \xi_{n,k} \quad \mathrm{for} \quad  t \in S_{0,0} =[0,1]\, .
\nonumber
\end{equation}
We need to study the path-wise convergence of the partial sums $X^{N}$ on $\Omega$ to see in which sense we can consider $\lim_{N \to \infty} X^{N}$ as a proper stochastic process. \\
Again, we only consider Gaussian Markov processes for which the functions $f$ and $g$ belong to the set of continuous functions $C[0,1]$.
If we designate the $L^{\infty}$ norms  of $f$ and $g$ on $[0,1]$ by ${\lVert f \rVert}_{\infty}$ and ${\lVert g \rVert}_{\infty}$, we can show that
\begin{eqnarray} \label{eq:UpBound}
\sup_{0 \leq k < 2^{n-1}} \sup_{ 0 \leq t \leq 1 }  \quad \big\vert \, \Psi_{n,k}(t)  \big\vert    \: & \leq  & \:   \sqrt{\frac{( r_{n,k} -m_{n,k}) ( m_{n,k} -l_{n,k})}{r_{n,k} -l_{n,k}}} \cdot {\lVert g \rVert}_{\infty}  \cdot {\lVert f \rVert}_{\infty}  \nonumber\\
& \leq & \hspace{40pt} 2^{-\frac{n+1}{2}} \cdot {\lVert g \rVert}_{\infty}  \cdot {\lVert f \rVert}_{\infty} \, . 
\end{eqnarray}
For $(f,g)$-bounded Gaussian Markov processes, this inequality provides us with the same upper bound to the elements $\Phi_{n,k}$ as in the case of a Wiener process times a constant ${\lVert g \rVert}_{\infty}  \cdot {\lVert f \rVert}_{\infty}$. 
By the same Borel-Cantelli argument as for the Haar construction of the Wiener process~\cite{Karatzas}, for almost every $\omega$ in $\Omega$, the sample path converges almost surely normally and uniformly in $t$ to a function $t \mapsto \overline{X}_{t}(\omega)$ when $N$ goes to infinity.
\\
It is worth noticing that, since $f$ and $g$ are continuous functions, so are the basis functions $\Psi_{n,k}$.
Then, for every $\omega$ in $\Omega$, the sample path $t \mapsto X_{t}^{N}(\omega)$ is a continuous function in $C[0,1]$.
As the convergence when $N$ tends to infinity is normal and uniform in t, the limit functions $t \mapsto \overline{X}_{t}(\omega)$ results to be in $C[0,1]$ almost surely on $\Omega$.
This allows us to define on $\Omega$ a limit process $\overline{X} = \lim_{N \to \infty} X^{N}$ with continuous paths.\\
Showing that $\overline{X}$ is an admissible discrete representation of the Gaussian Markov process $X$ only amounts to demonstrate that, for any integers $k$ and positive reals $t_{1} < t_{2}< \cdots <t_{k}$, the random vector $(\overline{X}_{t_{1}},\overline{X}_{t_{2}},\cdots,\overline{X}_{t_{k}})$ has a the same joint distribution as $(X_{t_{1}},X_{t_{2}},\cdots,X_{t_{k}})$.
As $\overline{X}$ is defined as the path-wise almost sure limit of $X^{N}$ when $N$ tends to infinity, this result is implied by the convergence in distribution of the continuous processes $X^{N}$ toward their limit $X$ ~\cite{Karatzas,Billingsley}.
We will therefore establish the convergence in distribution  of our representation in section \ref{sec:WeakConv}  and incidentally validate $\overline{X}$ as an exact representation of $X$.
In that perspective, we devote the following section to set out the meaning of the convergence in distribution  of our candidate process $X^{N}$.


\section{The Convergence in Distribution} \label{sec:ConvRep}


\subsection{The Finite-dimensional Measures}

In this section, we specify the finite-dimensional probability measures $P^{N}$ induced by the processes $X^{N}$. 
Beforehand, we introduce the notations 
\begin{equation}
[n,k] = 2^{n\!-\!1}\!+\!k \quad \mathrm{with} \quad 
\left\{
\begin{array}{ll}
 [0,0] = 0 \\
 {[0,1]} = 2^{N}
\end{array}
\right.
\nonumber
\end{equation}
to allow us to list the midpoints $m_{n,k}$ of the tree of supports $S_{n,k}$ in the prefix order. By reindexing according to
\begin{equation}
t_{[n,k]} = m_{n,k}=(2k\!+\!1) \, 2^{-n} \, ,
\nonumber
\end{equation}
we get an ordered sequence $t_{0}< \! t_{1} \! < \! \cdots \! < \! t_{2^{N}}$.
Let us now consider $C^{N}$ the finite-dimensional space of admissible functions for $X^{N}$ 
\begin{equation}
C^{N} = \mathrm{Vect}\left( {\Big\lbrace\Psi_{n,k}\Big\rbrace}_{ \substack{ \hspace{-8pt} \scriptscriptstyle 0 \leq n \leq N \\ \scriptscriptstyle 0\leq k<2^{n\!-\!1}}}\right) 
\nonumber
\end{equation}
When it is equipped with the $L^{\infty}$ norm, $C^{N}$ is a complete, separable metric space under the distance $d(f,g) = {\Vert f-g \Vert}_{\infty}$.
We can provide the space $C^{N}$ with the $\sigma$-algebra $\mathcal{B}(C^{N})$ generated by the cylinder sets $C_{B_{0},  \cdots,   B_{2^{N}}}$, which are defined for any collection of Borel sets $B_{0}, B_{1},  \cdots,  B_{2^{N}}$ in $\mathcal{B}(\mathbb{R})$ by
\begin{equation}
C_{B_{0},  \cdots,   B_{2^{N}}} = \Big \lbrace x \in C^{N} \; \vert \; x(t_{[n,k]}) \in B_{[n,k]}, \, 0 \leq n \leq N, \, 0 \leq k < 2^{n\!-\!1} \Big \rbrace \,
\nonumber
\end{equation}
The random process $X^{N}$ induces a natural  measure $P^{N}$ on $(C^{N},\mathcal{B}(C^{N}))$, such that
\begin{equation}
\forall B \in \mathcal{B}(C^{N}) \, , \quad P^{N}(B) = \mathrm{\bold{P}}(\omega \; \vert \; X^{N}(\omega) \in B)\, .
\nonumber
\end{equation}
Since for any $x$ in $C^{N}$ we have $x(0)=0$, the induced measure $P^{N}$ is entirely determined on the cylinder sets of the form $C_{B_{0},  \cdots,   B_{2^{N}}}$  with $B_{0}=\lbrace 0 \rbrace$.
Keeping this in mind, we show in annex that $P_N$ admits a probability density $p^{N}$:
for any cylinder set $C_{B_{0},  \cdots,   B_{2^{N}}}$ of $\mathcal{B}(C^{N})$ with $B_{0}=\lbrace 0 \rbrace$ we have
\begin{equation}
P^{N}(C_{B_{0},  \cdots,   B_{2^{N}}}) = \int_{B_{1}} \cdots \int_{B_{2^{N}}} p^{N}(x_{1}, \cdots ,  x_{2^{N}}) \,dx_{1}\cdots dx_{2^{N}} 
\nonumber
\end{equation}
where $p^{N}$  is made explicit with the help of the transition kernel $p$ of $X$ defined in \eqref{eq:CondProb}
\begin{equation}
p^{N}(x_{1}, \cdots ,  x_{2^{N}}) = \prod_{k=0}^{2^{N}\!-\!1} p(x_{k\!+\!1}, t_{k\!+\!1} \vert x_{k},t_{k})  \, .
\nonumber
\end{equation}
We want to specify in which sense the finite-dimensional probability measures $P^{N}$ defined on $(C^{N},\mathcal{B}(C^{N}))$ converge to a limit measure $P$ associated to $X$.


\subsection{The Weak Convergence}

Here, we consider that the stochastic processes $X^{N}$ and $X$ take value in the Wiener space, that is the space of continuous functions $C[0,1]$. 
This allows us to characterize the $X$-induced measure $\mu=P$ associated with $X$ on $C[0,1]$.
Defining the $X^{N}$-induced measures $\mu^{N}$ on $C[0,1]$ as well, we then state the convergence of $P^{N}$ toward $P$ in terms of weak convergence of $\mu^{N}$ toward $\mu$.\\
The process $X^{N}$ defined on some probability space  $\left( \Omega, \mathcal{F}, \mathrm{\bold{P}} \right)$ have continuous sample paths $t \mapsto X^{N}_{t}(\omega)$ and so does the Gaussian Markov process $X$.
Being a complete, separable metric space under the distance $d(f,g) = {\Vert f-g \Vert}_{\infty}$, the Wiener space is a natural space to define $X^{N}$- and $X$-induced measures.
We consider $X^{N}$ and $X$ as a random variables with values in the measurable space $\left( C[0,1], \mathcal{B}\left(C [0,1] \right) \right)$, where  $\mathcal{B}\left(C [0,1] \right)$ is the $\sigma$-field generated by  the cylinder sets of $C[0,1]$:
as previously $X^{N}$ and $X$ naturally induce the probability measures $\mu^{N}$ and $\mu$  defined by
\begin{equation} \label{eq:DefIndMeas}
\mu^{N}(B) = \mathrm{\bold{P}} \lbrace \omega \in \Omega \big \vert X^{N}(\omega) \in B \rbrace 
\; \: \: \mathrm{and} \; \: \:
\mu(B) = \mathrm{\bold{P}} \lbrace \omega \in \Omega \big \vert X(\omega) \in B \rbrace 
\nonumber
\end{equation} 
for any $B$ in $\mathcal{B}\left(C [0,1] \right)$.\\
More specially, the measure $\mu$ is called the Wiener measure of the Gaussian Markov process $X$.
Assuming a general cylinder set to be 
\begin{equation}
C^{t_{1}, \cdots, t_{n}}_{B_{1}, \cdots, B_{n}} = \big \lbrace x \in C[0,1] \; \big \vert \; x(t_{k}) \in B_{k}, \, 0 \! < \!  k \! \leq \! n \big \rbrace 
\nonumber
\end{equation}
for any $0 \! < \! t_{1} \! < \! t_{2} \! < \! \cdots \! < \! t_{n}\! \leq \! 1$ and any Borel sets $B_{1}, B_{2},\cdots, B_{n}$ in $\mathcal{B}(\mathbb{R})$, $\mu$ is the unique probability measure such that for any set $C^{t_{1}, \cdots, t_{n}}_{B_{1}, \cdots, B_{n}}$
\begin{eqnarray}
\mu(C^{t_{1}, \cdots, t_{n}}_{B_{1}, \cdots, B_{n}})
=
\int_{B_{1}} \! \! \! \cdots \! \int_{B_{n}} 
p(x_{1},t_{1} \vert 0, 0)
\cdots 
p(x_{n},t_{n} \vert x_{n\!-\!1}, t_{n\!-\!1}) 
\: dx_{t_{1}} \cdots dx_{t_{n}} \, .
\nonumber
\end{eqnarray}\\
Through the induced measures $\mu^{N}$ and $\mu$, we want to study the convergence of the process $X^{N}$ toward $X$ from a probabilistic point of view.
In that respect, the most general form of convergence one might expect is the weak convergence.
We recall that $\mu^{N}$ is weakly convergent to $\mu$ if and only if for any bounded continuous function of $C[0,1]$, we have
\begin{equation}
\lim_{N \to \infty} \int_{C} \phi(x) \, d\mu^{N}(x) = \int_{C} \phi(x) \, d\mu(x) \, ,
\nonumber
\end{equation}
where $C$ is a short notation for $C[0,1]$.


\subsection{The Convergence in Distribution}

We can now translate on the Wiener space  the much desirable property that the continuous processes $X^{N}$ converges in distribution toward the Gaussian Markov process $X$, which we denote
\begin{equation}
X^{N} \stackrel{\mathcal{D}}{\longrightarrow} X \, .
\nonumber
\end{equation}
We recall that, by definition, $X^{N}$ converges in distribution to $X$ if and only if for any bounded continuous function $\phi$ in $C[0,1]$ we have 
\begin{equation}
\lim_{N \to \infty} E_{N} \big ( \phi(X^{N}) \big ) =  E \big ( \phi(X) \big ) \, 
\nonumber
\end{equation}
where $E_{N}$ and $E$ are the expectations with respect to $P^{N}$ and $P$ respectively.
With the definitions made in the previous section, the convergence in distribution of the representation $X^{N}$ is rigorously equivalent to the weak convergence of the measures $\mu^{N}$ toward $\mu$ on the Wiener space $\left( C[0,1], \mathcal{B}\left(C [0,1] \right) \right)$.  \\
We will show this point in section \ref{sec:WeakConv} following the usual two steps reasoning inspired by the Prohorov theorem~\cite{Billingsley}.
To show the weak convergence of the sequence of measure $\mu^{N}$, it is enough to prove the two statements:
\begin{enumerate}
\item For every integer  $k>0$ and reals $0 \! \leq  \!  t_{1} \! <  \! t_{2} \!  < \! \cdots \! <  \! t_{k}  \! \leq  \! 1$, the finite-dimensional vector $(X^{N}_{t_{1}},X^{N}_{t_{2}},\cdots,X^{N}_{t_{k}})$ converges in distribution to $(X_{t_{1}},$ $X_{t_{2}},\cdots,X_{t_{k}})$ when $N$ tends to infinity.
\vspace{5pt}
\item The family of induced measures $\mu^{N}$  is tight: for every $\eta>0$, there exist a compact $K \subset C[0,1]$ such that $\mu^{N}(K) \geq 1-\eta$ for every $N$.
\end{enumerate}
Before establishing these two criteria, we first need to compute the limit of the covariance of $X^{N}$ when $N$ tends to infinity.
This calculation is the crucial point to validate our representation and we will carefully detail it in the following section.


\section{The Covariance Calculation} \label{sec:CovCal}


\subsection{Definition of the Auxiliary Basis }

Let us remember that the Gaussian Markov process $X$ admits a Doob's representation \eqref{eq:IntGM}.
If we posit appropriate regularity properties for $f$ and $g$, it is straigtforward to see that, by It\=o formula, such a process is solution of the stochastic differential equation
\begin{equation} \label{eq:Ito}
d \left( \frac{X_{t}}{g(t)} \right) = f(t) \cdot dW_{t} \, .
\nonumber
\end{equation}
It is then tempting to inject the proposed basis element $\Psi_{n,k}$ in the previous equation and to consider the functions $\Phi_{n,k}$ defined as
\begin{equation}
\Phi_{n,k}(t) = \frac{1}{f(t)} \frac{d}{du} {\left( \frac{\Psi_{n,k}(u)}{g(u)} \right)}_{u=t} \, .
\nonumber
\end{equation}
The functions $\Phi_{n,k}$ are actually well-defined despite the division by $f$, which is potentially zero, since calculations show that the $\Phi_{n,k}$ are given explicitly for $n>0$ by
\begin{equation} \label{eq:Phi definition}
\Phi_{n,k}(t)= \left\{ \begin{array}{lllll}
\displaystyle &&L_{n,k} \cdot f(t)
\quad \textrm{if \: $l_{n,k}  \leq  t  <  m_{n,k}$} \: , \\
\\
\displaystyle &-& R_{n,k}  \cdot f(t)
\quad \textrm{if \: $m_{n,k}  \leq t < r_{n,k}$} \: ,\\
\\
 &&0  \quad \textrm{otherwise} \: . \end{array} \right.
\end{equation}
As for the element $\Psi_{0,0}$, it gives rise to the well-defined function
\begin{equation}
\Phi_{0,0}(t) = \frac{f(t)}{\sqrt{h(r_{0,0})-h(l_{0,0})}} \,.
\nonumber
\end{equation}
We will show that the family of functions $\Phi_{n,k}$ is a Hilbert system of a subspace of $L^{2}(0,1)$, a property that will prove useful to compute, for any $t$ and $s$ in $[0,1]$, the limit of the covariance  $E(X_{t}^{N} \cdot X_{s}^{N})$ when $N$ tends to infinity.
The consideration of this result will also enable us to interpret the coefficients $L_{n,k}$ and $R_{n,k}$. \\


\subsection{Characterization as an Hilbert System}

From now on, we denote $(f,g)$ the usual inner product of $f$ and $g$ in $L^2(0,1)$.
Let us introduce $\overline{{\mathcal{E}}_{f}}$, the closure of the vectorial space ${\mathcal{E}}_{f}$  defined as  
\begin{equation}
{\mathcal{E}}_{f} = \lbrace \phi \in L^2(0,1)\, \vert \, \exists \, \varphi \in L^2(0,1), \phi=f\cdot\varphi\rbrace \, .
\nonumber
\end{equation}
With the usual inner product, the space $\overline{{\mathcal{E}}_{f}}$ inherits the structure of a Hilbert space from $L^2(0,1)$ for being a closed subspace of a Hilbert space.
It is immediate to see that the $\Phi_{n,k}$ belong to ${\mathcal{E}}_{f}$ when written as \eqref{eq:Phi definition}.
We need to prove that the $\Phi_{n,k}$ constitute an orthonormal family for the usual inner product and that the vectorial space of their finite linear combinations  is dense in  $\overline{{\mathcal{E}}_{f}}$.\\
Let us start with the orthonormal property of the family and consider two functions $\Phi_{n,k}$ and $\Phi_{n',k'}$ with $(n,k) \ne (n',k')$.
If $n=n'$ and $k \ne k'$, $\Phi_{n,k}$ and $\Phi_{n',k'}$ have disjoint supports and their inner product is necessarily zero. 
Assuming that $n' > n$, $\Phi_{n,k}$ and $\Phi_{n',k'}$ have intersecting supports $S_{n,k}$ and $S_{n',k'}$ if and only if $S_{n,k}$ is strictly included in $S_{n',k'}$.
Then, it is very useful to remark the nullity of the following inner product
\begin{eqnarray}
 \left( f,\Phi_{n,k} \right)
 & = &
 L_{n,k} \cdot \int_{l_{n,k}}^{m_{n,k}} f^{2}(u) \,du+ R_{n,k} \cdot \int_{m_{n,k}}^{r_{n,k}} f^{2}(u) \, du 
 \nonumber\\
 & = & 
\sqrt{
\frac{ \big( h(r_{n,k})-h(l_{n,k}) \big) 
  \big( h(m_{n,k})-h(l_{n,k}) \big)}
 { h(r_{n,k})-h(m_{n,k})}
 }
 - 
 \nonumber\\
&  &  - 
\sqrt{
\frac{ \big( h(r_{n,k})-h(l_{n,k}) \big) 
  \big( h(m_{n,k})-h(l_{n,k}) \big)}
 { h(r_{n,k})-h(m_{n,k})}
 } = 0 \, ,
 \nonumber
\end{eqnarray}
which entails that $\Phi_{n,k}$ and $\Phi_{n',k'}$ are orthogonal for $n \ne n'$.
As for the norm of the functions $\Phi_{n,k}$, we directly compute for $n>0$
\begin{eqnarray}
\left( \Phi_{n,k},\Phi_{n,k} \right)
 & = &
 L^{2}_{n,k} \cdot \int_{l_{n,k}}^{m_{n,k}} f^{2}(u) \, du+ R^{2}_{n,k} \cdot \int_{m_{n,k}}^{r_{n,k}} f^{2}(u) \, du 
 \nonumber\\
 & = & 
 \frac{h(r_{n,k})-h(m_{n,k})}{h(r_{n,k})-h(l_{n,k})} + \frac{h(m_{n,k})-h(l_{n,k})}{h(r_{n,k})-h(l_{n,k})} = 1 \, 
 \nonumber
\end{eqnarray}
and it is straightforward to see that $\left( \Phi_{0,0},\Phi_{0,0} \right) = 1$ for $n=0$.
Hence, we have proved that the collection of $\Phi_{n,k}$ forms an orthonormal family of functions in ${\mathcal{E}}_{f}$.\\
We still have to show that the linear combinations of $\Phi_{n,k}$ generate a dense vectorial space in $\overline{{\mathcal{E}}_{f}}$.
We can easily be convinced of this point once we consider the family of functions $f \cdot H_{n,k}$, where the $H_{n,k}$ designate the Haar functions adapted to the supports  $S_{n,k}$.
The orthonormal system $H_{n,k}$ is dense in $L^{2}(0,1)$ as soon as condition \eqref{eq:DenseCond} is satisfied.
Then, as each $f \cdot H_{n,k}$ can be obtained by finite linear combination of $\Phi_{n,k}$, we conclude that 
\begin{equation}
\overline{\mathrm{Vect}\left( {\Big\lbrace\Phi_{n,k}\Big\rbrace}_{ \substack{ \hspace{-17pt} \scriptscriptstyle n\geq0 \\ \scriptscriptstyle 0\leq k<2^{n\!-\!1}}}\right)} = \overline{{\mathcal{E}}_{f}} \, ,
\nonumber
\end{equation}
The interpretation of the coefficients $R_{n,k}$ and $L_{n,k}$ is made conspicuous in this context. The natural surjective morphism of $L^{2}[0,1]$ to $\overline{{\mathcal{E}}_{f}}$ is the application $ \phi \mapsto f \cdot \phi$. The family $f \cdot H_{n,k}$ in  $\overline{{\mathcal{E}}_{f}}$ is the image of the Haar basis $H_{n,k}$, but it is not a Hilbert system of $\overline{{\mathcal{E}}_{f}}$.
The coefficients $R_{n,k}$ and $L_{n,k}$ results from the operations of orthonormalization of the family $f \cdot H_{n,k}$ to form a Hilbert system of $\overline{{\mathcal{E}}_{f}}$.


\subsection{Application of the Parseval Relation}

Now that these preliminary remarks have been made, we can evaluate, for any $t$ and $s$ in $[0,1]$, the limit of the covariance $E(X_{t}^{N} \cdot X_{s}^{N})$ when $N$ tends to infinity.
As $\xi_{n,k}$ are independent Gaussian random variables of normal law $\mathcal{N}(0,1)$, we see that the covariance of $X^{N}$ is given by
\begin{eqnarray} \label{eq: VarPart}
E \left( X_{t}^{N} \cdot X_{s}^{N} \right) = 
\sum_{n=0}^{N} \sum_{ \hspace{5pt} 0 \leq k < 2^{n\!-\!1} } \; \Psi_{n,k}(t) \cdot \Psi_{n,k}(s) \, .
\end{eqnarray}
To compute the limit of the right-hand side in \eqref{eq: VarPart}, we need to remark that the element of the basis $\Psi_{n,k}$ and the function of the auxiliary Hilbert system $\Phi_{n,k}$ are linked by the following relation
\begin{eqnarray}\label{eq:IntPsi}
\Psi_{n,k}(t) & = & g(t) \cdot \int_{0}^{1} \chi_{\lbrack 0, t \rbrack}(u) \frac{d}{dv} {\left( \frac{\Psi_{n,k}(v)}{g(v)} \right)}_{v=u} \, du
\nonumber\\
& = &   \int_{0}^{1} \chi_{\lbrack 0, t \rbrack}(u) \, g(t) f(u) \cdot \Phi_{n,k}(u)\, du \, .
\end{eqnarray}
In the previous expression, we use the indicator functions of the segment $\lbrack 0,t\rbrack$ defined as
\begin{equation} \label{eq:indicator}
\chi_{ \left[ 0, t \right] }(u)= \left\{ \begin{array}{ll}
1& \textrm{if $0 \leq u \leq t $}\\
0  & \textrm{otherwise} \quad \end{array} \right. \: .
\nonumber
\end{equation}
In order to simplify the notations, we now introduce the functions $\eta_{t}$ elements of the Hilbert space ${\mathcal{E}}_{f}$
\begin{equation}
\eta_{t}(u) = \chi_{\lbrack 0, t \rbrack}(u) \, g(t) f(u) \, .
\nonumber
\end{equation}
With the help of the function $\eta_{t}$, we can then write the integral definition of the basis element $\Psi_{n,k}$ \eqref{eq:IntPsi} as the inner product in $\overline{{\mathcal{E}}_{f}}$
\begin{equation}
\Psi_{n,k}(t) = \left( \eta_{t}, \Phi_{n,k}\right) \, .
\nonumber
\end{equation}
Remembering that the family of functions $\Phi_{n,k}$ is a Hilbert system of $\overline{{\mathcal{E}}_{f}}$, we can make use of Parseval identity, which reads
\begin{eqnarray} \label{eq:Parseval}
\left( \eta_{s}, \eta_{t}\right) =
\sum_{n=0}^{\infty} \sum_{ \hspace{5pt} 0 \leq k < 2^{n\!-\!1} } 
\left( \eta_{s}, \Phi_{n,k}\right) 
 \cdot 
\left( \eta_{t}, \Phi_{n,k}\right)  \, .
\end{eqnarray}
Thanks to this relation, we can conclude the evaluation of the variance of $\overline{X}$, since a direct explicitation of \eqref{eq:Parseval} yields  
\begin{eqnarray} \label{eq:ParsExpl}
\int_{0}^{1} \chi_{\lbrack 0, t \rbrack}(u) \, g(t) f(u) \cdot \chi_{\lbrack 0, s \rbrack}(u) \, g(s) f(u)\, du
=
\sum_{n=0}^{\infty} \sum_{ \hspace{5pt} 0 \leq k < 2^{n\!-\!1} } \; \Psi_{n,k}(t) \cdot \Psi_{n,k}(s) \, . \nonumber
\\
\end{eqnarray}
The left term in \eqref{eq:ParsExpl} precisely happens to be the same as the covariance of $X$ given by relation \eqref{eq:XCov} and we recap the statement by saying
\begin{displaymath}
\lim_{N \to \infty} E \left( X^{N}_{t} \cdot X^{N}_{s} \right)  =    g(s)g(t) \cdot h\big(\min(t,s) \big)  =   E \left( X_{t} \cdot X_{s} \right) \, .
\end{displaymath}
In section \ref{sec:WeakConv}, we will use this relation  to show that the random vector $(X^{N}_{t_{1}},X^{N}_{t_{2}},$ $\cdots,X^{N}_{t_{k}})$ converges in distribution toward $(X_{t_{1}},$ $X_{t_{2}},\cdots,X_{t_{k}})$ for any reals $0 \! \leq  \!  t_{1} \! <  \! t_{2} \!  < \! \cdots \! <  \! t_{k}  \! \leq  \! 1$.


\section{The Convergence of the Representation}\label{sec:WeakConv}


\subsection{Formulation of the two Sufficient Criteria}

We are now in a position to proceed to the demonstration of the weak convergence of the induced measures $\mu^{N}$ toward $\mu$.
This will prove the convergence in distribution of the representation
\begin{equation}
X^{N} \stackrel{\mathcal{D}}{\longrightarrow} X \, ,
\nonumber
\end{equation}
and validate the fact that the limit process $\overline{X}$ is an exact representation of $X$.
We just need to verify that our representation satisfies the two criteria mentioned in section \ref{sec:ConvRep}.\\
First, we show that the family of measures $\mu_{N}$ fulfill the tightness condition, using a characterization of tightness on $C[0,1]$.
This characterization is a version of the Arzel\`a-Ascoli theorem.
It is formulated in term of $mod(x,\delta)$ the $\delta$-modulus of continuity of a function $x$ in $C[0,1]$ :
\begin{equation}
\forall \delta>0, \quad mod(x,\delta) = \max_{\vert t-s \vert \leq \delta} \big \vert x(t)-x(s)  \big \vert \, .
\nonumber
\end{equation}
Bearing in mind that we only consider the subset of functions $x$ in $C[0,1]$ for which $x(0) \! = \! 0$, the characterization reads 
\begin{equation} \label{eq:TightCrit}
\lim_{\delta \to 0} \; \sup_{N \in \mathbb{N}} \; \mu^{N}\Big( x \in C[0,1] \; \big \vert \; \vert mod\big(x, \delta \big) \vert > \epsilon \Big) = 0 \, .
\end{equation}\\
Then, for any reals $0 \! \leq  \!  t_{1} \! <  \! t_{2} \!  < \! \cdots \! <  \! t_{k}  \! \leq  \! 1$, we show the convergence in distribution of the random vectors $(X^{N}_{t_{1}},X^{N}_{t_{2}},\cdots,X^{N}_{t_{r}})$
\begin{equation}
\big( X^{N}_{t_{1}},X^{N}_{t_{2}}, \cdots , X^{N}_{t_{r}}\big) \stackrel{\mathcal{D}}{\longrightarrow} \big( X_{t_{1}},X_{t_{2}}, \cdots , X_{t_{r}}\big) \, .
\nonumber
\end{equation}
Dealing with finite-dimensional vectors, we can use the Cram\'er-Wold device~\cite{Karatzas}: 
if we designate $\mathscr{C}^{N}$ and $\mathscr{C}$ the characteristic functions of $(X^{N}_{t_{1}},X^{N}_{t_{2}},$ $\cdots,X^{N}_{t_{r}})$ and $(X_{t_{1}},X_{t_{2}},\cdots,X_{t_{r}})$ defined on $\mathbb{R}^{r}$, we just have to show the point-wise convergence of the characteristic functions  $\mathscr{C}^{N}$ toward $\mathscr{C}$.


\subsection{The Tightness of the Induced Family of Distributions}

We will confirm that the family of induced measure $\mu^{N}$ satisfy the tightness criterion \eqref{eq:TightCrit} on $C[0,1]$.
First, define the random variables
\begin{equation}
b_{n} = \max_{0 \leq k <2^{n-1}} \vert \xi_{n,k} \vert \, .
\nonumber
\end{equation}
We apply the usual Borel-Cantelli lemma and state:
there exists a set $\Omega'$ with $\mathrm{\bold{P}}(\Omega')=1$ such that for every $\omega$ in $\Omega'$, there is an integer $n(\omega)$ such that for all $n>n(\omega)$ we have $b_{n}(\omega) \leq n$.
Let us set
\begin{equation}
\Omega_{n}'=\lbrace \omega \in \Omega' \; \vert \; n(\omega) \leq n \rbrace
\nonumber
\end{equation}
It clearly defines an increasing sequence of sets $\Omega_{0}' \subset \Omega_{1}' \subset \cdots \subset \Omega_{n}'$ with $\lim_{n \to \infty}\mathrm{\bold{P}}(\Omega_{n}')=1$.
For any $\eta > 0$, there is a $n(\eta)$ such that for all $n>n(\eta)$, we have $\mathrm{\bold{P}}(\Omega_{n}') > 1-\eta/2$.
Then Considering the existence of the upper bound \eqref{eq:UpBound}, for every $\omega$ in $\Omega_{n}'$ with $N>n(\eta)$, we can write
\begin{eqnarray} \label{eq:Tight1}
\lefteqn{ \hspace{-20pt}
\big \vert X^{N}_{t}(\omega) - X^{N}_{s}(\omega) \big \vert 
\quad \leq \quad
\sum_{n=0}^{n(\eta)} \sum_{ \hspace{5pt} 0 \leq k < 2^{n(\eta)\!-\!1} } \big \vert \Psi_{n,k}(t) - \Psi_{n,k}(s) \big \vert \big \vert \xi_{n,k}(\omega) \big \vert }
 \nonumber\\
 && \hspace{100pt}
+
\sum_{n=n(\eta)+1}^{N} 2n \cdot 2^{-\frac{n+1}{2}} \cdot {\lVert g \rVert}_{\infty}  \cdot {\lVert f \rVert}_{\infty}
\end{eqnarray}
The previous inequality actually holds for any $N>0$ with the conventions of setting the elements $\Psi_{n,k}$ to zero when $n>N$.
For all $\epsilon > 0$ and all $\eta > 0$ there exists an integer $n(\epsilon, \eta)>n(\eta)$ such that 
\begin{equation}
\sum_{n=n(\epsilon, \eta)+1}^{\infty} 2n \cdot 2^{-\frac{n+1}{2}} \cdot {\lVert g \rVert}_{\infty}  \cdot {\lVert f \rVert}_{\infty}
\leq
\frac{\epsilon}{2}
\nonumber
\end{equation}
Now, for any given $\epsilon > 0$ and $\eta>0$, we chose a real $M(\epsilon,\eta)>0$ large enough so that for any $\omega$ in $\Omega$
\begin{equation}
\mathrm{\bold{P}} \bigg(\max_{n \leq n(\epsilon, \eta)}b_{n}(\omega) >  M(\epsilon,\eta) \bigg)
\quad \leq \quad
\frac{\eta}{2}
\nonumber
\end{equation}
and we finally define the set $\Omega_{\epsilon,\eta}$ as
\begin{equation}
\Omega_{\epsilon,\eta} = \bigg \lbrace \omega \in \Omega'_{n(\epsilon,\eta)} \; \big \vert \;  \max_{n \leq n(\epsilon, \eta)}b_{n}(\omega) \leq M(\epsilon,\eta) \bigg \rbrace \, .
\nonumber
\end{equation}
Every element $\Psi_{n,k}$ is a continuous function defined on a compat support $S_{n,k}$ in $[0,1]$.
As a result, the finite set of functions $\Psi_{n,k}$ for $n\leq n(\epsilon,\eta)$ is uniformly equicontinuous:
for any given $\epsilon > 0$ and $\eta>0$, there is  a real $\delta(\epsilon,\eta)>0$ such that 
\begin{eqnarray}
\lefteqn{
\forall t,s \in [0,1], \quad \vert t -s\vert \leq \delta(\epsilon,\eta) \quad \Rightarrow  \quad 
}
\nonumber\\
&& 
\forall n,k \quad \mathrm{with} \:
\left\{
\begin{array}{ll}
 0 \leq n\leq n(\epsilon,\eta) \\
  0 \leq k < 2^{n-1}
\end{array}
\right. ,
\,\big \vert \Psi_{n,k}(t) - \Psi_{n,k}(s) \big \vert \leq \frac{\epsilon}{2^{n(\epsilon,\eta)\!+\!1}M(\epsilon,\eta)} \, .
\nonumber
\end{eqnarray}
For all $\epsilon > 0$ and all $\eta > 0$, writing the inequality \eqref{eq:Tight1} on $\Omega_{\epsilon,\eta}$ yields to
\begin{equation}
\forall \omega \in \Omega_{\epsilon,\eta}, \quad \delta <  \delta(\epsilon,\eta) \quad \Rightarrow  \quad \forall N \in \mathbb{N}, \quad \big \vert mod\big(X^{N}(\omega), \delta \big) \big \vert \leq \epsilon 
\nonumber
\end{equation}
As we have defined $\Omega_{\epsilon,\eta}$ so that $\mathrm{\bold{P}}(\Omega_{\epsilon,\eta}) \geq 1-\eta$, we have shown that 
\begin{equation}
\forall \epsilon, \eta >0, \quad  \delta <  \delta(\epsilon,\eta) \quad \Rightarrow \quad \sup_{N \in \mathbb{N}} \mathrm{\bold{P}} \bigg( \big \vert mod\big(X^{N}(\omega), \delta \big) \big \vert > \epsilon \bigg)
\leq 
\eta \, ,
\nonumber
\end{equation}
which proves the tightness criterion on the Wiener space C[0,1] for the family of continuous processes $X^{N}$.\\


\subsection{Point-wise Convergence of the Characteristic Functions}

For any $( \lambda_{1},\lambda_{2}, \cdots , \lambda_{r})$ in $\mathbb{R}^{r}$, evaluating $\mathscr{C}^{N}$, the characteristic function of $(X^{N}_{t_{1}},X^{N}_{t_{2}},$ $\cdots,X^{N}_{t_{r}})$, yields 
\begin{eqnarray}
 \mathscr{C}^{N}( \lambda_{1},\lambda_{2}, \cdots , \lambda_{r})
&=&
\int_{C} e^{i \sum_{p=1}^{r} \lambda_{p} x(t_{p}) } \, d\mu^{N}(x) \, ,
\nonumber
\end{eqnarray}
We use the definition \eqref{eq:SimpSum} of the partial sum $X^{N}$ in terms of the basis elements $\Psi_{n,k}$ to explicit the calculation of $ \mathscr{C}^{N}$ on $C^{N}$, the space of admissible functions for $X^{N}$:
\begin{eqnarray}
 \mathscr{C}^{N}( \lambda_{1},\lambda_{2}, \cdots , \lambda_{r})
&=&
\int_{C^{N}} e^{i \sum_{p=1}^{r} \lambda_{p} X^{N}_{t_{p}} } \, \mathrm{\bold{P}}(dX^{N}) \, .
\nonumber
\end{eqnarray}
Remember that each coefficient $\xi_{n,k}$ is independently distributed according to a normal law $\mathcal{N}(0,1)$ in the representation of $X^{N}$. We have then
\begin{equation}
 \mathscr{C}^{N}( \lambda_{1},\lambda_{2}, \cdots , \lambda_{r})
= 
\quad \prod_{n=0}^{N} \prod_{\hspace{5pt}  0 \leq k < 2^{n-1}} \int_{\mathbb{R}} \: e^{i \xi_{n,k} \left(\sum_{p}^{r} \lambda_{p} \Psi_{n,k}(t_{p}) \right)} \, \mathrm{\bold{P}}(d\xi_{n,k}) \, .
\nonumber
\end{equation}
Therefore, we can compute each terms of the previous product:
\begin{eqnarray}
\int_{\mathbb{R}} \: e^{i \xi_{n,k} ( \lambda_{p} \Psi_{n,k}(t_{p}) )} \, \mathrm{\bold{P}}(d\xi_{n,k})
&=&
\int_{\mathbb{R}} \: e^{i \xi_{n,k} ( \lambda_{p} \Psi_{n,k}(t_{p}) )} \, \frac{1}{\sqrt{2\pi}}e^{-\frac{\xi^{2}_{n,k}}{2}} \, d\xi_{n,k}
\nonumber\\
&=&
\exp{\Big(- \frac{1}{2} {\big(\lambda_{p} \Psi_{n,k}(t_{p})\big)}^{2} \Big)} \, .
\nonumber
\end{eqnarray}
Back to the formulation of $\mathscr{C}_{N}$, we end up with the analytical expression
\begin{equation}
 \mathscr{C}^{N}( \lambda_{1},\lambda_{2}, \cdots , \lambda_{r})
 = 
\exp{\Big(- \frac{1}{2} \sum_{p}^{r} \sum_{q}^{r} \: \lambda_{p} \lambda_{q} \: \rho^{N}(t_{p},s_{q}) \Big)} \, ,
\nonumber
\end{equation}
where we have used the short notation $\rho^{N}$ for the covariance of $X^{N}$
\begin{equation}
\rho^{N}(t,s)
=
E( X^{N}_{t} \cdot X^{N}_{s})
=
\sum_{n=0}^{N} \sum_{\hspace{5pt}  0 \leq k < 2^{n-1}}   \Psi_{n,k}(t) \cdot \Psi_{n,k}(s) \, .
\nonumber
\end{equation}
With the covariance calculation of section \ref{sec:CovCal}, we have demonstrated that for  $t$, $s$ in $[0,1]$
\begin{equation}
\lim_{N \to \infty} \rho^{N}(t,s) 
=
E(X_{t} \cdot X_{s})
=
g(t)g(s)\cdot h\big(\min(t,s)\big)
=
\rho(t,s) \, .
\nonumber
\end{equation}
For any $( \lambda_{1},\lambda_{2}, \cdots , \lambda_{r})$ in $\mathbb{R}^{r}$, it directly entails the point-wise convergence of $\mathscr{C}^{N}$
\begin{eqnarray}
\lim_{N \to \infty}   \mathscr{C}^{N}( \lambda_{1},\lambda_{2}, \cdots , \lambda_{r})
&=& 
\exp{\Big(- \frac{1}{2} \sum_{p}^{r} \sum_{q}^{r} \: \lambda_{p} \lambda_{q} \: \rho(t_{p},s_{q}) \Big)} 
\nonumber\\
&\stackrel{def}{=}& 
 \mathscr{C}( \lambda_{1},\lambda_{2}, \cdots , \lambda_{r}) \, ,
\nonumber
\end{eqnarray}
where we remark that $\mathscr{C}$ is the characteristic function of the random vector $( X_{t_{1}},X_{t_{2}}, \cdots , X_{t_{k}})$.\\
By the Cram\'er-Wold device, the point-wise convergence of $\mathscr{C}^{N}$ toward $\mathscr{C}$ proves the convergence in distribution of any finite-dimensional random vectors $X^{N}_{t_{1}},X^{N}_{t_{2}}, \cdots , X^{N}_{t_{r}}$ toward the random vector $X_{t_{1}},X_{t_{2}}, \cdots , X_{t_{r}}$
\begin{equation}
\big( X^{N}_{t_{1}},X^{N}_{t_{2}}, \cdots , X^{N}_{t_{r}}\big) \stackrel{\mathcal{D}}{\longrightarrow} \big( X_{t_{1}},X_{t_{2}}, \cdots , X_{t_{r}}\big) \, .
\nonumber
\end{equation}
It also proves that $\overline{X}$ is an exact representation of the Gaussian Markov process $X$ since, as the almost-sure path-wise limit of $X^{N}$ when $N$ tends to infinity, $\overline{X}$ follows the same law as the law of the process $X$.
As we have already established the tightness of the family of induced distributions $\mu^{N}$, it is enough to prove the weak convergence of $\mu^{N}$ toward $\mu$ on the Wiener space, which is equivalent to the convergence in distribution of the continuous process $X^{N}$ towards $X$.


\appendix

\begin{scriptsize}


\section{Gaussian calculation}\label{appA}

We want to establish the analytical results given in section \ref{AnalyticalResults}.
Assuming $t_{x} < t_{y}< t_{z}$, we want to compute $p( X_{t_{y}}   =  y \, \vert X_{t_{x}}  =  x, X_{t_{z}}  =  z)$. We use the Markov property as mentioned in the third section

\begin{eqnarray} 
 p( X_{t_{y}}   =  y \, \vert X_{t_{x}}  =  x, X_{t_{z}}  =  z)  & = &
\frac{ p( X_{t_{y}} =  y\,  \vert X_{t_{x}}  = x) 
\cdot 
p( X_{t_{z}}  =  z \, \vert X_{t_{y}}  = y) }
{  p( X_{t_{z}}  =  z \, \vert X_{t_{x}} =  x ) } 
\nonumber
\end{eqnarray} 

\noindent to get its analytical expression in terms of $g$ and $h$

\begin{eqnarray} \label{eq1}
\frac{
\displaystyle
\frac{1}{g(t_y)\sqrt{2 \pi} \cdot  \sigma_{x,y} } 
\exp{ \left( - \frac{ \left( \frac{y}{g(t_y)} - \frac{x}{g(t_x)} \right)^{2} }{2 \cdot \sigma_{x,y}^2} \right) }
\cdot
\frac{1}{g(t_z)\sqrt{2 \pi} \cdot  \sigma_{y,z} } 
\exp{  \left( - \frac{ \left( \frac{z}{g(t_z)} - \frac{y}{g(t_y)}  \right)^{2} }{2 \cdot \sigma_{y,z}^2} \right) }
}
{
\displaystyle
\frac{1}{g(t_z)\sqrt{2 \pi} \cdot  \sigma_{x,z} } 
\exp{ \left( - \frac{ \left( \frac{z}{g(t_z)} - \frac{x}{g(t_x)}  \right)^{2} }{2 \cdot \sigma_{x,z}^2} \right) }
}
\nonumber\\
\: ,
\end{eqnarray}

\noindent with the expressions of the variance of $X_t$ between any two consecutive points

\begin{equation}
\sigma_{x,y}^2 = h(t_y)-h(t_x) \: , \quad \sigma_{y,z}^2 = h(t_z)-h(t_y) \: , \quad \sigma_{x,z}^2 = h(t_z)-h(t_x) \: .
\end{equation}

\noindent After factorization of the exponentials, the resulting exponent of expression \eqref{eq1} is written

\begin{equation}
\begin{split}
-& 
\frac{1}{2 \cdot \sigma_{x,y}^2} \left( \frac{y}{g(t_y)} - \frac{x}{g(t_x)} \right)^{2} 
-
\frac{1}{2 \cdot \sigma_{y,z}^2} \left( \frac{z}{g(t_z)} - \frac{y}{g(t_y)} \right)^{2} 
+
\frac{1}{2 \cdot \sigma_{x,z}^2} \left( \frac{z}{g(t_z)} - \frac{x}{g(t_x)} \right)^{2} 
= 
\\
\\
&-
\underbrace{
\frac{1}{2 \cdot g^{2}(t_y)} \left( \frac{1}{\sigma_{x,y}^2} + \frac{1}{\sigma_{y,z}^2}  \right) 
}_{C_y}  \; y^2
+
\frac{1}{g(t_y)} \left(  \frac{ x }{ g(t_x)\sigma_{x,y}^2} + \frac{ z }{ g(t_z)\sigma_{y,z}^2} \right)\; y
\\
& \hspace{120pt} -
\frac{ x^{2}}{2 \cdot g^{2}(t_x) \sigma_{x,y}^{2}}
-
\frac{ z^{2}}{2 \cdot g^{2}(t_z) \sigma_{y,z}^{2}}
+
\frac{ \left( \frac{z}{g(t_z)} - \frac{x}{g(t_x)} \right)^{2} }{2 \cdot \sigma_{x,z}^2} 
\nonumber \: .
\end{split}
\end{equation}

\noindent We then factorize the term $C_y$ so that we can write the exponent in the form

\begin{eqnarray}
\lefteqn{
- 
\frac{1}{2 \cdot g^{2}(t_y)} \cdot \frac{\sigma^{2}_{x,z}}{\sigma^{2}_{x,y} \cdot \sigma^{2}_{y,z}}
\left\{
y^{2} 
- 
2 \left( \frac{g(t_y)}{g(t_x)} \cdot \frac{\sigma_{y,z}^{2} }{\sigma_{x,z}^{2}} \cdot x + \frac{g(t_y)}{g(t_z)} \cdot \frac{\sigma_{x,y}^{2}}{\sigma_{x,z}^{2}} \cdot z \right)
y
\right\} 
}
\nonumber\\
&&
+
\frac{1}{2 \cdot \sigma_{x,z}^2} \left( \frac{z}{g(t_z)} - \frac{x}{g(t_x)} \right)^{2}  
-
\frac{ x^{2}}{2 \cdot g^{2}(t_x) \sigma_{x,y}^{2}}
-
\frac{ z^{2}}{2 \cdot g^{2}(t_z) \sigma_{y,z}^{2}}
\nonumber
\end{eqnarray}

\noindent to obtain the canonical expression

\begin{multline}
- 
\frac{1}{2 \cdot g^{2}(t_y)} \cdot \frac{\sigma^{2}_{x,z}}{\sigma^{2}_{x,y} \cdot \sigma^{2}_{y,z}}
\left[
y
-
\left(  \frac{g(t_y)}{g(t_x)} \cdot \frac{\sigma_{y,z}^{2} }{\sigma_{x,z}^{2}} \cdot x + \frac{g(t_y)}{g(t_z)} \cdot \frac{\sigma_{x,y}^{2}}{\sigma_{x,z}^{2}} \cdot z \right)
\right]^{2}
+ 
\frac{1}{2 \cdot \sigma_{x,z}^2} \left( \frac{z}{g(t_z)} - \frac{x}{g(t_x)} \right)^{2} 
\\
+ \:
\underbrace{
\frac{1}{2 \cdot g^{2}(t_y)} \cdot \frac{\sigma^{2}_{x,z}}{\sigma^{2}_{x,y}\sigma^{2}_{y,z}}
\left(  \frac{g(t_y)}{g(t_x)} \cdot \frac{\sigma_{y,z}^{2} }{\sigma_{x,z}^{2}} \cdot x + \frac{g(t_y)}{g(t_z)} \cdot \frac{\sigma_{x,y}^{2}}{\sigma_{x,z}^{2}} \cdot z \right)^{2} 
-
\frac{ x^{2}}{2 \cdot g^{2}(t_x) \sigma_{x,y}^{2}}
-
\frac{ z^{2}}{2 \cdot g^{2}(t_z) \sigma_{y,z}^{2}}
}_{Q}
\nonumber \: .
\end{multline}

\noindent The quantity $Q$ can be further simplify after being expanded

\begin{eqnarray}
\begin{split}
Q  
& = 
\frac{1}{2\cdot g^{2}(t_x)} 
\left(
\frac{\sigma^{2}_{y,z}}{\sigma^{2}_{x,y}\cdot \sigma^{2}_{x,z}} 
-
\frac{1}{\sigma_{x,y}^{2}}
\right)
x^2
+
\frac{1}{g(t_x)g(t_z)} \cdot \frac{1}{\sigma^{2}_{x,z}} \; \; xz
\\
&\quad +
\frac{1}{2\cdot g^{2}(t_z)} 
\left(
\frac{\sigma^{2}_{x,y}}{\sigma^{2}_{y,z}\cdot \sigma^{2}_{x,z}} \cdot z^2
-
\frac{1}{\sigma_{y,z}^{2}}
\right)
z^2
\\
& = 
-\frac{1}{2\cdot g^{2}(t_x)} \cdot \frac{\sigma^{2}_{x,y}}{\sigma^{2}_{x,z}} \; \; x^2
+
\frac{1}{g(t_x)g(t_z)} \cdot \frac{1}{\sigma^{2}_{x,z}} \; \; xz
-
\frac{1}{2\cdot g^{2}(t_z)} \cdot \frac{\sigma^{2}_{y,z}}{\sigma^{2}_{x,z}} \; \; z^2
\\
& = 
- \frac{1}{2 \cdot \sigma_{x,z}^{2}} {\left( \frac{z}{g(t_z)} - \frac{x}{g(t_x)} \right)}^{2}
\nonumber
\end{split}
\end{eqnarray}

\noindent All terms except quadratic ones in $y$ cancel out in the exponent, giving

\begin{eqnarray} \label{eq2}
\lefteqn{
p( X_{t_{y}}   =  y \, \vert X_{t_{x}}  =  x, X_{t_{z}}  =  z)   = 
}
\nonumber\\
&&
\frac{1}{\sqrt{2 \pi} \cdot \: \scriptstyle g(t_{y}) \textstyle \frac{ \sigma_{x,y} \sigma_{y,z} }{\sigma_{x,z}} } 
\cdot
\exp{
\left(
-
\frac{
\left[
y
-
\left(  \frac{g(t_y)}{g(t_x)} \cdot \frac{\sigma_{y,z}^{2} }{\sigma_{x,z}^{2}} \cdot x + \frac{g(t_y)}{g(t_z)} \cdot \frac{\sigma_{x,y}^{2}}{\sigma_{x,z}^{2}} \cdot z \right)
\right]^{2}
}
{
2  \cdot \: \scriptstyle g^{2}(t_y)\textstyle \frac{\sigma_{x,y}^{2} \cdot \sigma_{y,z}^{2}}{\sigma_{x,z}^{2}}
}
\right)
} \: .
\end{eqnarray}

\noindent We  sum up the expression \eqref{eq2} noticing it represents the distribution of a normal law $\mathcal{N}({}_{ \scriptscriptstyle X \displaystyle} \mu(t_{y}),{}_{ \scriptscriptstyle X \displaystyle} \sigma(t_{y}))$, whose parameters read

\begin{equation} 
{}_{ \scriptscriptstyle X \displaystyle} \mu(t_{y}) 
= 
\frac{g(t_y)}{g(t_x)} \cdot \frac{\sigma_{y,z}^{2} }{\sigma_{x,z}^{2}} \cdot x + \frac{g(t_y)}{g(t_z)} \cdot \frac{\sigma_{x,y}^{2}}{\sigma_{x,z}^{2}} \cdot z 
= 
\frac{g(t_y)}{g(t_x)} \cdot \frac{ h(t_{z})-h(t_{y}) }{ h(t_{z})-h(t_{x}) } \cdot x + \frac{g(t_y)}{g(t_z)} \cdot \frac{ h(t_{y})-h(t_{x}) }{ h(t_{z})-h(t_{x}) } \cdot z 
\nonumber
\end{equation}

\begin{equation} 
{}_{ \scriptscriptstyle X \displaystyle} \sigma(t_{y})^{2}
= 
g^{2}(t_{y}) \cdot \frac{\sigma_{x,y}^{2} \cdot \sigma_{y,z}^{2}}{\sigma_{x,z}^{2}}
=
g^{2}(t_{y}) \cdot \frac{\big( h(t_{y})-h(t_{x}) \big) \big( h(t_{y})-h(t_{x}) \big)}{h(t_{z})-h(t_{x})}
\nonumber
\end{equation}


\section{Induced measures}\label{appB}


\subsection*{\scriptsize  Finite-dimensional measures}

We first recall the definition of the $2^{N}$-dimensional vectorial space $C^{N}$ defined as
\begin{equation}
C^{N} = \mathrm{Vect}\left( {\Big\lbrace\Psi_{n,k}\Big\rbrace}_{ \substack{ \hspace{-8pt} \scriptscriptstyle 0 \leq n \leq N \\ \scriptscriptstyle 0\leq k<2^{n\!-\!1}}}\right) \, .
\nonumber
\end{equation}

\noindent We specify a given element $x$ of $C^{N}$ by writting

\begin{equation}
x(t) =  \sum_{n=0}^{N} \sum_{ \hspace{5pt} 0 \leq k < 2^{n\!-\!1} } \Psi_{n,k}(t) \cdot a_{n,k} \quad \mathrm{with} \quad a_{n,k} \in \mathbb{R}\, ,
\nonumber
 \end{equation}

\noindent and we remark that $x$ can be viewed as a sample path of the process $X^{N}$ if we posit $\xi_{n,k}(\omega) = a_{n,k} $.
We then introduce the following notations

\begin{equation}
[n,k] = 2^{n\!-\!1}\!+\!k \quad \mathrm{with} \quad 
\left\{
\begin{array}{ll}
 [0,0] = 0 \\
 {[0,1]} = 2^{N}
\end{array}
\right.
\nonumber
\end{equation}

\noindent to enumerate indices between $0$ and $2^{N}$ in the prefix order.
With this convention, we introduce the positive reals $0 = t_{0} < t_{1}< \cdots < t_{2^{N}} = 1$ defined as

 \begin{equation}
t_{[n,k]} = m_{n,k}=(2k+1) \, 2^{-n} \quad \mathrm{for} \quad 0 < n < N \quad \mathrm{and} \quad 0 \leq k < n \, .
\nonumber
\end{equation}

\noindent For any collection of Borel sets $B_{0}, B_{1}, \cdots,  B_{2^{N}}$ in $\mathcal{B}(\mathbb{R})$, we define the cylinder sets

\begin{equation}
C_{B_{0},  \cdots,   B_{2^{N}}} = \lbrace x \in C^{N} \vert x(t_{k}) \in B_{k} , \, 0< k \leq 2^{N}\rbrace
\nonumber
\end{equation}

\noindent which generate the $\sigma$-algebra $\mathcal{B}(C^{N})$.\\
We want to show that there exists a probability density $p^{N}$ on $\mathbb{R}^{2^{N}}$ such that we have for any sets $\lbrace 0 \rbrace = B_{0}, B_{1}, \cdots,  B_{2^{N}}$ in $\mathcal{B}(\mathbb{R})$

\begin{equation}
\mathrm{\bold{P}}(x \in C_{B_{0},B_{1},  \cdots,   B_{2^{N}}}) = \int_{B_{1}} \cdots \int_{B_{2^{N}}} p^{N}(x_{1}, \cdots ,  x_{2^{N}}) \,dx_{1}\cdots dx_{2^{N}} 
\nonumber
\end{equation}

\noindent and that $p^{N}$ has the following simple expression  in terms of the transition kernel $p$ associated with the Gaussian Markov process $X$

\begin{equation} \label{eq:probabilityDensity}
p^{N}(x_{1}, \cdots ,  x_{2^{N}}) = \prod_{k=0}^{2^{N}\!-\!1} p\Big( (x_{k\!+\!1}, t_{k\!+\!1}) \Big \vert (x_{k},t_{k}) \Big) \, .
\end{equation}

\noindent We start noticing that, by construction of $X^{N}$, we have

\begin{equation}
\mathrm{\bold{P}}(\xi_{n,k} \in da_{n,k}) 
= 
\mathrm{\bold{P}}\Big( X_{m_{n,k}} - \mu(m_{n,k}) \in dx_{[n,k]}  \Big \vert X_{l_{n,k}} = x_{[n\!-\!1,2k]} , X_{r_{n,k}} = x_{[n\!-\!1,2k\!+\!1]} \Big) 
\nonumber
\end{equation}

\noindent with the usual definitions for $l_{n,k}$, $r_{n,k}$ and $m_{n,k}$ and with $ dx_{[n,k]} = \sigma_{n,k} \cdot da_{n,k} $.
 As the Gaussian variables $\xi_{n,k}$ are all independent in the definition of $X^{N}$ and remembering the definition of the probability density $q$ in section \ref{AnalyticalResults}, we define $p^{N}$ for $N \geq 0$ as

\begin{eqnarray}
\lefteqn{
p^{N}(x_{1}, \cdots ,  x_{2^{N}}) = p\Big((x_{2^{N}},r_{1,0}) \Big \vert (x_{0},l_{1,0})  \Big) q\Big((x_{[0,1]} ,m_{1,0}) \Big \vert (x_{0},l_{1,0}), (x_{2^{N}},r_{1,0}) \Big)
\cdots}
\nonumber\\
&&
\hspace{80pt}
\cdots \prod_{k=0}^{2^{N\!-\!1}-1} q\Big((x_{[N,k]} ,m_{N,k}) \Big \vert (x_{[N\!-\!1,2k]},l_{N,k}), (x_{[N\!-\!1,2k\!+\!1]},r_{N,k}) \Big)
\nonumber \, ,
\end{eqnarray}

\noindent and we will show by recurrence on $N$ that $p^{N}$ as the expected expression \eqref{eq:probabilityDensity}.\\
The basis statement is obvious for $N = 0$.
As for the inductive step, let us write the probability density $p^{N\!+\!1}$ as

\begin{eqnarray}
\lefteqn{
p^{N\!+\!1}(x_{1}, \cdots ,  x_{2^{N\!+\!1}}) 
=}
\nonumber\\
&&
p^{N}(x_{1}, \cdots ,  x_{2^{N}})
\prod_{k=0}^{2^{N}\!-\!1} 
q\Big((x_{[N\!+\!1,k]} ,m_{N\!+\!1,k}) \Big \vert (x_{[N,2k]},l_{N\!+\!1,k}), (x_{[N,2k\!+\!1]},r_{N\!+\!1,k}) \Big) \, .
\nonumber
\end{eqnarray}

\noindent By our recurrence hypothesis, we have then

\begin{eqnarray}
\lefteqn{
p^{N\!+\!1}(x_{1}, \cdots ,  x_{2^{N\!+\!1}})
=}
\nonumber\\
&&
\prod_{k=0}^{2^{N\!-\!1}} p\Big( (x_{k\!+\!1}, t_{k\!+\!1}) \Big \vert (x_{k},t_{k}) \Big)
\prod_{k=0}^{2^{N}\!-\!1} 
q\Big((x_{[N\!+\!1,k]} ,m_{N\!+\!1,k}) \Big \vert (x_{[N,2k]},l_{N\!+\!1,k}), (x_{[N,2k\!+\!1]},r_{N\!+\!1,k}) \Big) \, .
\nonumber
\end{eqnarray}

\noindent But by definition of $q$, we have 

\begin{eqnarray}
\lefteqn{
q\Big((x_{[N\!+\!1,k]} ,m_{N\!+\!1,k}) \Big \vert (x_{[N,2k]},l_{N\!+\!1,k}), (x_{[N,2k\!+\!1]},r_{N\!+\!1,k}) \Big)
=}
\nonumber\\
&&
\frac{
p\Big( (x_{[N,2k\!+\!1]},r_{N\!+\!1,k}) \Big \vert (x_{[N\!+\!1,k]} ,m_{N\!+\!1,k}) \Big)
p\Big( (x_{[N\!+\!1,k]} ,m_{N\!+\!1,k}) \Big \vert (x_{[N,2k]},l_{N\!+\!1,k}) \Big)
}
{
p\Big( (x_{[N,2k\!+\!1]},r_{N\!+\!1,k}) \Big \vert (x_{[N,2k]},l_{N\!+\!1,k})) \Big)
} \, .
\nonumber
\end{eqnarray}

\noindent The product of the denominators in the previous expressions cancels out the term $p^{N}$ in $p^{N\!+\!1}$ so that we have proven that

\begin{equation}
p^{N\!+\!1}( x_{1}, \cdots ,  x_{2^{N\!+\!1}}) 
= 
\prod_{k=0}^{2^{N\!+\!1}\!-\!1} p\Big( (x_{k\!+\!1}, t_{k\!+\!1}) \Big \vert (x_{k},t_{k}) \Big) \, .
\nonumber
\end{equation}


\subsection*{\scriptsize Characteristic Functionals}

We want to express the characteristic functionals of the induced measures  $\mu^{N}$ and of the Wiener measure $\mu$.
In view of this, we recall that, by the Riesz representation theorem, the dual space of $C[0,1]$ is the space $M[0,1]$ of finite measures on $[0,1]$.
For $x$ in $C[0,1]$ and $\nu$ in $M[0,1]$, we write the duality product

\begin{equation}
\langle x , \, \nu \rangle = \int_{0}^{1} x(t) \, d\nu(t) \, .
\nonumber
\end{equation}

\noindent Considering first the induced measure $\mu_{N}$, the characteristic functional $\mathscr{C}_{N}$ is defined on $M[0,1]$ as the Fourier transform of $\mu_{N}$ on $C[0,1]$

\begin{equation}
\mathscr{C}_{N}(\nu) =  \int_{C} e^{i \langle x , \, \nu \rangle} \, d\mu^{N}(x) \quad \mathrm{with} \quad \nu \in M[0,1] \, .
\nonumber
\end{equation}

\noindent By construction of the process $X^{N}$ and independence of the random variables $\xi_{n,k}$, we have

\begin{equation}
\mathscr{C}_{N}(\nu) =  \int_{C^{N}} e^{i \langle X^{N} \!\! , \, \nu \rangle} \, \mathrm{\bold{P}}(dX^{N}) = \prod_{n=0}^{N} \prod_{\hspace{5pt}  0 \leq k < 2^{n-1}} \int_{\mathbb{R}} \: e^{i \xi_{n,k} \langle \Psi_{n,k} , \, \nu \rangle} \, \mathrm{\bold{P}}(d\xi_{n,k}) \, .
\nonumber
\end{equation}

\noindent  The random variables $\xi_{n,k}$ being Gaussian of law $\mathcal{N}(0,1)$, we furthermore have

\begin{equation}
\int_{\mathbb{R}} \: e^{i \xi_{n,k} \langle \Psi_{n,k} , \, \nu \rangle} \, \mathrm{\bold{P}}(d\xi_{n,k})
=
\int_{\mathbb{R}} \: e^{i \xi_{n,k} \langle \Psi_{n,k} , \, \nu \rangle} \, \frac{1}{\sqrt{2\pi}}e^{-\frac{\xi^{2}_{n,k}}{2}} \, d\xi_{n,k}
=
e^{- \frac{1}{2} {\langle \Psi_{n,k} , \, \nu \rangle}^{2}} \, .
\nonumber
\end{equation}

\noindent This allows to finally write the expression of the characteristic functional $\mathscr{C}_{N}$

\begin{eqnarray}
\mathscr{C}_{N}(\nu) 
=
\exp{\bigg(- \frac{1}{2} \sum_{n=0}^{N} \sum_{\hspace{5pt}  0 \leq k < 2^{n-1}} {\big \langle \Psi_{n,k} , \, \nu \big \rangle}^{2} \bigg)} \, .
\nonumber
\end{eqnarray}

\noindent If we express the product of duality, we can formulate the functional $\mathscr{C}_{N}$ in its common form~\cite{Donsker}

\begin{eqnarray} \label{eq:CharAnalytic}
\mathscr{C}_{N}(\nu) 
&=& 
\exp{\bigg(- \frac{1}{2} \sum_{n=0}^{N} \sum_{\hspace{5pt}  0 \leq k < 2^{n-1}}  \int_{0}^{1}  \Psi_{n,k}(t) \, d\nu(t)  \int_{0}^{1}  \Psi_{n,k}(s) \, d\nu(s) \bigg)}
\nonumber\\
&=&
\exp{\bigg(- \frac{1}{2} \int_{0}^{1} \int_{0}^{1} \rho^{N}(t,s) \,  d\nu(t)d\nu(s) \bigg)}
\nonumber
\end{eqnarray}

\noindent where we made apparent $\rho^{N}$, the correlation function of the process $X^{N}$:

\begin{equation}
\rho^{N}(t,s)
=
E( X^{N}_{t} \cdot X^{N}_{s})
=
\sum_{n=0}^{N} \sum_{\hspace{5pt}  0 \leq k < 2^{n-1}}   \Psi_{n,k}(t) \cdot \Psi_{n,k}(s) \, .
\nonumber
\end{equation}

\noindent From the convergence results of our expansion, we directly have that the characteristic functional $\mathscr{C}$ of the Wiener measure $\nu$

\begin{equation} \label{eq:CharactFunction}
\mathscr{C}(\nu)
=
\exp{\bigg(- \frac{1}{2} \int_{0}^{1} \int_{0}^{1} \rho(t,s) \,  d\nu(t)d\nu(s) \bigg)} \, ,
\end{equation}

\noindent  with $\rho$ the continuous correlation function of the Gaussian Markov process $X$:

\begin{equation}
\rho(t,s)
=
E(X_{t} \cdot X_{s})
=
 g(t)g(s)\cdot h\big(\min(t,s)\big) \, .
 \nonumber
\end{equation}

\noindent  Now let $t_{1}<t_{2}<\cdots<t_{k}$ be some reals in $[0,1]$, we posit the measure $\nu$ as

\begin{equation}
\nu = \sum_{p=1}^{k} \lambda_{p} \delta_{t_{p}} \, ,
\nonumber
\end{equation}

\noindent  where $ \delta_{t_{p}}$ denotes the Dirac distribution concentrated at $t_{p}$.
If we inject the expression of $\nu$ in the result \eqref{eq:CharactFunction}, we find the expression of the characteristic function of  $X_{t_{1}},X_{t_{2}}, \cdots , X_{t_{k}}$ as one might expect.

\noindent Then, assume the distribution $\nu$ admits a density $\theta$ in $L^{2}(0,1)$ with respect to the Lebesgues measure, we have

\begin{equation}
\big \langle \Psi_{n,k} ,  \nu \big \rangle = \int_{0}^{1} \Psi_{n,k}(t) \theta(t) \, dt = \big( \Psi_{n,k} , \, \theta \big) 
\nonumber
\end{equation}

\noindent Thanks to the auxiliary orthonormal basis $\Phi_{n,k}$, we can further write

\begin{equation}
\big( \Psi_{n,k} , \, \theta \big)  
=
\int_{0}^{1}
\underbrace{\frac{1}{f(t)} \frac{d}{du} {\left( \frac{\Psi_{n,k}(u)}{g(u)} \right)}_{u=t} }_{\Phi_{n,k}}
 \cdot
 \:
 f(t) \int_{0}^{t} g(u) \theta(u) \, du 
 \quad
 dt \: ,
 \nonumber
\end{equation}

\noindent which directly leads to the following simple expression for the characteristic functional

\begin{equation} 
\mathscr{C}(\nu)
=
\exp{\bigg(- \frac{1}{2} \int_{0}^{1}  {\Big( f(t) \int_{0}^{t} g(s) \theta(s) \, ds \Big)}^{2} \,  dt \bigg)} \, .
\nonumber
\end{equation}

\end{scriptsize}


\bibliography{AAP-Paper}

\end{document}